\documentclass[12pt,a4paper]{article}
\usepackage{amsmath,amsfonts,amssymb,amscd}

\def\al{\alpha}

\def\rk{\operatorname{rk}}
\def\H{\operatorname{H}}
\def\gr{\operatorname{gr}}
\def\Ker{\operatorname{Ker}}
\def\Ad{\operatorname{Ad}}
\def\ad{\operatorname{ad}}

\def\id{\operatorname{id}}
\def\ev{\operatorname{ev}}
\def\red{\operatorname{red}}
\def\d{\operatorname{d}}
\def\GL{\operatorname{GL}}
\def\OSp{\operatorname{OSp}}
\def\Sp{\operatorname{Sp}}
\def\Q{\operatorname{Q}}
\def\pr{\operatorname{pr}}
\def\Im{\operatorname{Im}}
\def\pt{\operatorname{pt}}
\def\Hom{\operatorname{Hom}}

\def\Lie{\operatorname{Lie}}

\newcounter{th}
\def\t{\refstepcounter{th}{\bf \noindent{Theorem} \arabic{th}. }}

\newcounter{le}
\def\l{\refstepcounter{le}{\bf \noindent{Proposition} \arabic{le}. }}

\newcounter{lem}
\def\lem{\refstepcounter{lem}{\bf \noindent{Lemma} \arabic{lem}. }}

\newcounter{de}

\newcounter{ex}
\def\ex{\refstepcounter{ex}{\bf \noindent{Example} \arabic{ex}. }}

\begin{document}

\begin{center}
{\Large {\bf On holomorphic functions on a compact complex
homogeneous
supermanifold} \footnote{Work supported by SFB $\mid$ TR 12. 
}}
\end{center}

\begin{center}
     E.G. Vishnyakova
\end{center}

\noindent\textsc{Abstract.} It is well-known that non-constant
holomorphic functions do not exist on a compact complex manifold.
This statement is false for a supermanifold with a compact
reduction.  In this paper we study the question under what
conditions non-constant holomorphic functions do not exist on a
 compact homogeneous complex supermanifold. We describe also the
vector bundles determined by split homogeneous complex
supermanifolds.

As an application, we compute the algebra of holomorphic functions
on the classical flag supermanifolds which were introduced in
\cite{Man}.

\bigskip

\begin{center}
{\bf 1. Preliminaries}
\end{center}

\begin{center}
{\it 1.1 Lie supergroups and homogeneous supermanifolds}
\end{center}

We will use the word "supermanifold" \,in the sense of Berezin and
Leites (see \cite{BL, ley}).  All the time, we will be interested in
the complex-analytic version of the theory. Let $(M,\mathcal{O}_M)$
be a supermanifold. 
The underlying complex manifold $M$ is called the {\it reduction} of
$(M,\mathcal{O}_M)$. The superalgebra $H^0(\mathcal{O}_M)$ is called
the {\it superalgebra of (global) holomorphic functions} on
$(M,\mathcal{O}_M)$. A function $f\in H^0(\mathcal{O}_M)$ is called
{\it constant} if $f|_U$ does not depend on even and odd coordinates
for every coordinate superdomain $(U,\mathcal{O}_M|_U)\subset
(M,\mathcal{O}_M)$.

\smallskip

\ex\label{example1} Let $\mathcal{E}$ be a locally free sheaf on
$M$. Then $(M,\bigwedge \mathcal{E})$ is a supermanifold. Let
$U\subset M$ be a coordinate domain of $M$ with coordinates $(x_i)$.
Assume that $\mathcal{E}|_U$ is free and $(\xi_j)$ is a local basis.
Then $(U,\bigwedge \mathcal{E}|_{U})$ is a superdomain with
coordinates $(x_i,\xi_j)$. Note that any $f\in H^0(\bigwedge^p
\mathcal{E})\backslash \{0\}$, where $p>0$, is not constant. Suppose
that $M$ is compact. Obviously, it does not follow that
$H^0(\bigwedge^p \mathcal{E})= \{0\}$ for $p>0$.

\smallskip

We denote by $\mathcal{J}_M\subset\mathcal{O}_M$ the subsheaf of
ideals generated by odd elements of the structure sheaf. The sheaf
$\mathcal{O}_M/\mathcal{J}_M$ is naturally identified with the
structure sheaf $\mathcal{F}_M$ of $M$. The natural homomorphism
$\mathcal{O}_M\to\mathcal{F}_M$ will be denoted by $f\mapsto
f_{\red}$.
 A morphism
$\phi:(M,\mathcal{O}_M)\to (N,\mathcal{O}_N)$ of supermanifolds will
be denoted by $\phi = (\phi_{\red},\phi^*)$, where $\phi_{\red}:
M\to N$ is the corresponding mapping of the reductions and
$\phi^*:\mathcal{O}_N\to(\phi_{\red})_*(\mathcal{O}_M)$ is the
homomorphism of the structure sheaves. If $x\in M$ and $\mathfrak
m_x$ is the maximal ideal of the local superalgebra $(\mathcal
O_M)_x$, then the vector superspace $T_x(M,\mathcal O_M)=(\mathfrak
m_x/\mathfrak m_x^2)^*$ is the tangent space to $(M,\mathcal O_M)$
at $x\in M$. Denote by $\mathcal{T}_M$ the sheaf of derivations of
the structure sheaf $\mathcal{O}_M$. It is a sheaf of Lie
superalgebras with the Lie bracket $[X,Y]:=X\circ Y -
(-1)^{p(X)p(Y)}Y\circ X$, where $p(Z)$ is the parity of $Z$. We will
use the following notation
$\mathfrak{v}(M,\mathcal{O}_M)=H^0(\mathcal{T}_M)$ for the Lie
superalgebra of vector fields on $(M,\mathcal{O}_M)$. From the
inclusions $v(\mathfrak m_x)\subset (\mathcal O_M)_x$ and
$v(\mathfrak m_x^2)\subset \mathfrak m_x$, where $v\in \mathfrak
v(M,\mathcal O_M)$, it follows that $v$ induces an even linear
mapping $\ev_x(v):\mathfrak m_x/\mathfrak m_x^2\to (\mathcal
O_M)_x/\mathfrak m_x\simeq \Bbb C$. In other words, $\ev_x(v)\in
T_x(M,\mathcal O_M)$, and so we obtain an even linear map
\begin{equation}\label{ev}
\ev_x:\mathfrak v(M,\mathcal O_M)\to T_x(M,\mathcal O_M).
\end{equation}


 A {\it Lie supergroup} is a group object in
the category of supermanifolds, i.e., a supermanifold $(G,\mathcal
O_G)$, for which the following three morphisms are defined:
$\mu:(G,\mathcal O_G)\times (G,\mathcal O_G)\to (G,\mathcal O_G)$
(the multiplication morphism), $\iota:(G,\mathcal O_G)\to
(G,\mathcal O_G)$ (the inversion morphism),
$\varepsilon:(\pt,\mathbb C)\to (G,\mathcal O_G)$ (the identity
morphism). Moreover, these morphisms should satisfy the usual
conditions, modeling the group axioms. The underlying manifold $G$
 is a complex Lie group. The element
$e=\varepsilon_{\red}(\pt)$ is the identity element of $G$. We will
denote by $\mathfrak{g}$ the Lie superalgebra of $(G,\mathcal O_G)$.
By definition, $\mathfrak{g}$ is the subalgebra of
$\mathfrak{v}(G,\mathcal O_G)$ consisting of all right invariant
vector fields on $(G,\mathcal{O}_G)$. It is well known that any
right invariant vector field $Y$ has the form
\begin{equation}
\label{left inv vect field}
 Y=(X\otimes \id)\circ \mu^*
\end{equation}
 for a certain $X\in T_e(G,\mathcal{O}_G)$ and the map $X\mapsto (X\otimes \id)\circ \mu^*$ is an
isomorphism of the vector space $T_e(G,\mathcal{O}_G)$ onto
$\mathfrak{g}$, see \cite{Var}, Theorem $7.1.1$.

 An {\it action of a Lie supergroup
$(G,\mathcal O_G)$ on a supermanifold} $(M,\mathcal O_M)$ is a
morphism $\nu:(G,\mathcal O_G)\times (M,\mathcal O_M)\to (M,\mathcal
O_M)$ such that the following conditions hold:
\begin{itemize}
  \item $\nu \circ (\mu\times \id)=\nu\circ (\id\times \nu)$;
  \item $\nu\circ (\varepsilon\times \id)=\id$.
\end{itemize}
In this case $\nu_{\red}$ is the action of $G$ on $M$.

Let $\nu:(G,\mathcal O_G)\times (M,\mathcal O_M)\to (M,\mathcal
O_M)$ be an action. Then there is a homomorphism of the Lie
superalgebras $\overline{\nu}:\mathfrak g\to \mathfrak v(M,\mathcal
O_M)$ given by the formula
\begin{equation}\label{hom}
 X\mapsto
(X\otimes\, \id)\circ \nu^*.
\end{equation}

As in \cite{onipi}, we use the following definition of a transitive
action.
 An action $\nu$ is called {\it
transitive} if $\nu_{\red}$ is transitive and the mapping $\ev_{x}
\circ \overline{\nu}$ is surjective for all $x\in M$. (The map
$\ev_{x}$ is given by (\ref{ev}).) In this case the supermanifold
$(M,\mathcal O_M)$ is called $(G,\mathcal O_G)$-{\it homogeneous}.
 A supermanifold $(M,\mathcal O_M)$
is called {\it homogeneous} if it possesses a transitive action of a
certain Lie supergroup.

Suppose that a closed Lie subsupergroup $(H,\mathcal{O}_H)$ of
$(G,\mathcal{O}_G)$ (this means that the Lie subgroup $H$ is closed
in $G$) is given. Denote by $j$ the inclusion of $(H,\mathcal{O}_H)$
into $(G,\mathcal{O}_G)$. Consider the corresponding coset
superspace $(G/H,\mathcal{O}_{G/H})$, see \cite{Fio_Varad, Kostant}.
Denote by $\mu_{G\times H}$ the composition of the morphisms
$$
(G,\mathcal{O}_G)\times (H,\mathcal{O}_H)\stackrel{\id\times
j}{\longrightarrow} (G,\mathcal{O}_G)\times
(G,\mathcal{O}_G)\xrightarrow{\mu} (G,\mathcal{O}_G),
$$
 by $\pr_1:(G,\mathcal{O}_G)\times (H,\mathcal{O}_H)\to
(G,\mathcal{O}_G)$ the projection onto the first factor, and by
$\pi$ the natural mapping $G\to G/H$, $g\mapsto gH$. Let us take
$U\subset G/H$ open.  Then
\begin{equation}\label{sheaf O_G/H}
\mathcal{O}_{G/H}(U)= \{ f\in \mathcal{O}_G(\pi^{-1}(U))\mid
(\mu_{G\times H})^*(f)=\pr^*_1(f)\}.
\end{equation}
Denote by $\nu: (G,\mathcal{O}_G)\times (G/H,\mathcal{O}_{G/H})\to
(G/H,\mathcal{O}_{G/H})$ the natural action. It is given by
$\nu^*(f)=\mu^*(f)$, where $f\in \mathcal{O}_{G/H}(U)$. Hence if
$X\in \mathfrak{g}$, $f\in \mathcal{O}_{G/H}$, we have
$\overline{\nu}(X)(f)=X(f)$. Sometimes we will denote the
supermanifold $(G/H,\mathcal{O}_{G/H})$ also by
$(G,\mathcal{O}_G)/(H,\mathcal{O}_H)$.

\medskip

\ex\label{example2} Let $(G,\mathcal{O}_G)$ be a Lie supergroup, $H$
a Lie subgroup of $G$. Then $(H,\mathcal{F}_H)$ is also a Lie
subsupergroup of $(G,\mathcal{O}_G)$. It is well known that the
sheaf $\mathcal{O}_{G/H}$ is isomorphic to $\mathcal{F}_{G/H}
\otimes \bigwedge(\mathfrak{g}^*_{\bar 1})$, where
$\mathfrak{g}=\Lie (G,\mathcal O_G)$, see, e.g., \cite[Proposition
$2$]{V}. If $G/H$ is compact and connected,
$H^0(\mathcal{O}_{G/H})\simeq \bigwedge(\mathfrak{g}^*_{\bar 1})$.
Hence there are compact homogeneous complex supermanifolds with
non-constant holomorphic functions.

\medskip

\begin{center}
{\it 1.2 The Harish-Chandra pairs}
\end{center}

The structure sheaf of a Lie supergroup and the supergroup morphisms
can be explicitly described in terms of the corresponding Lie
superalgebra using so-called (super) Harish-Chandra pairs, see
\cite{Bern}. A {\it Harish-Chandra pair} is a pair
$(G,\mathfrak{g})$ that consists of a Lie group $G$ and a Lie
superalgebra $\mathfrak{g}=\mathfrak{g}_{\bar
0}\oplus\mathfrak{g}_{\bar 1}$, where $\mathfrak{g}_{\bar 0}$ is the
Lie algebra of $G$, provided with a representation $\al_G$ of $G$ in
$\mathfrak{g}$ such that
\begin{itemize}
  \item $\al_G$ preserves the parity and induces the adjoint representation
of $G$ in $\mathfrak{g}_{\bar 0}$,
  \item the differential $(\d\al_G)_e$ at the identity $e\in G$ coincides with
the adjoint representation $\ad$ of $\frak g_{\bar 0}$ in $\frak g$.
\end{itemize}

 Super Harish-Chandra pairs form a
category. (The definition of a morphism see in \cite{Bern}.)
The following theorem was proved in \cite{Kostant}.

\medskip
\t\label{main} {\it The category of real Lie supergroups is
equivalent to the category of real Harish-Chandra pairs.  }

\medskip

In the complex case, the equivalence of the categories was shown in
\cite{V}.

If a Harish-Chandra pair $(G,\mathfrak{g})$ is given, it determines
the Lie supergroup $(G,\widehat{\mathcal{O}}_G)$ in the following
way, see \cite{kosz}. Let $\mathfrak{U}(\mathfrak{g})$ be the
universal enveloping superalgebra of $\mathfrak{g}$. It is clear
that $\mathfrak{U}(\mathfrak{g})$ is a
$\mathfrak{U}(\mathfrak{g}_{\bar 0})$-module, where
$\mathfrak{U}(\mathfrak{g}_{\bar 0})$ is the universal enveloping
algebra of $\mathfrak{g}_{\bar 0}$. The natural action of
$\mathfrak{g}_{\bar 0}$ on the sheaf $\mathcal{F}_G$ gives rise to a
structure of $\mathfrak{U}(\mathfrak{g}_{\bar 0})$-module on
$\mathcal{F}_G(U)$ for any open set $U\subset G$. Putting
$$
\widehat{\mathcal{O}}_G(U) = \Hom_{\mathfrak{U}(\mathfrak{g}_{\bar
0})}(\mathfrak{U}(\mathfrak{g}), \mathcal{F}_G(U))
$$
for every open $U\subset G$, we get a sheaf
$\widehat{\mathcal{O}}_G$ of $\mathbb{Z}_2$-graded vector spaces
(here we assume that the functions from $\mathcal{F}_G(U)$ are
even). The enveloping superalgebra $\mathfrak{U}(\mathfrak{g})$ has
a Hopf superalgebra structure (see \cite{Scheu}). Using this
structure we can define the product of elements from
$\widehat{\mathcal{O}}_G$ such that $\widehat{\mathcal{O}}_G$
becomes a sheaf of superlgebras. A supermanifold structure on
$\widehat{\mathcal{O}}_G$ is determined by the isomorphism
$\widehat{\mathcal{O}}_G \stackrel{\sim}{\to} \Hom(\bigwedge
(\mathfrak{g}_{\bar 1}), \mathcal{F}_G)$, $f\mapsto f\circ \gamma$,
where
\begin{equation}
\label{isomorphism}
 \gamma: \bigwedge(\mathfrak{g}_{\bar 1})\to
\mathfrak{U}(\mathfrak{g}),\quad X_1\wedge \cdots \wedge X_r\mapsto
\frac{1}{r!}\sum_{\sigma\in S_r}(-1)^{|\sigma|} X_{\sigma(1)}\cdots
X_{\sigma(r)}.
\end{equation}
 The following formulas define the
multiplication morphism, the inversion morphism and the identity
morphism respectively (see \cite{BagSta}):
\begin{equation}\label{umnozh}
\begin{split}
\mu^*(f)(X\otimes Y)(g,h)&=f(X\cdot \al_G(g)(Y))(gh);\\
\iota^*(f)(X)(g)&=f(\al_G(g^{-1})(S(X)))(g^{-1});\\
\varepsilon^*(f)&=f(1)(e).
\end{split}
\end{equation}
Here $X,Y\in\mathfrak{U}(\mathfrak{g}),\, f\in
\widehat{\mathcal{O}}_G,\, g,\,h\in G$ and $S$ is the antipode map
of $\mathfrak{U}(\mathfrak{g})$. Here we identify the enveloping
superalgebra $\mathfrak{U}(\mathfrak{g}\oplus\mathfrak{g})$ with the
tensor product $\mathfrak{U}(\mathfrak{g})\otimes
\mathfrak{U}(\mathfrak{g})$.

A Harish-Chandra pair $(H,\mathfrak{h})$ is called a {\it
Harish-Chandra subpair} of a Harish-Chandra pair $(G,\mathfrak{g})$
if $H$ is a Lie subgroup of $G$ and $\mathfrak{h}$ is a  Lie
subsuperalgebra of $\mathfrak{g}$, s.t. $\mathfrak{h}_{\bar 0}=\Lie
H$ and $\alpha_H=\alpha_G|H$. There is a correspondence between
Harish-Chandra subpairs of $(G,\mathfrak{g})$ and Lie subsupergroups
of $(G,\mathcal O_G)$, see, e.g., \cite{V}. (The Lie supergroup
$(G,\mathcal O_G)$ corresponds to the Harish-Chandra pair
$(G,\mathfrak{g})$.)

Let $\nu$ be an action of $(G,\mathcal{O}_G)$ on
$(M,\mathcal{O}_M)$, $x\in M$ and $\delta_x:(\pt,\mathbb{C})\to
(M,\mathcal{O}_M)$ the morphism such that $\delta_x(\pt)=x$. Denote
by $\nu_x$ the following composition:
$$
(G,\mathcal{O}_G)\times (\pt,\mathbb{C})\stackrel{\id \times
\delta_x}{\longrightarrow} (G,\mathcal{O}_G)\times (M,\mathcal{O}_M)
\stackrel{\nu}{\to} (M,\mathcal{O}_M).
$$
Consider the Harish-Chandra subpair $(G_x,\mathfrak{g}_x)$ of
$(G,\mathfrak{g})$, $\mathfrak{g}=\Lie (G,\mathcal O_G)$, where
$G_x\subset G$ is the stabilizer of $x$ and $\mathfrak{g}_x=\Ker
(\d\nu_x)_e$. A subsupergroup $(G_x,\mathcal{O}_{G_x})$ is called
the {\it stabilizer} of $x$ if it is determined by
$(G_x,\mathfrak{g}_x)$.

Denote by $l_g$, $g\in G$, the following composition:
$$
(G,\mathcal{O}_G)=(g,\mathbb{C})\times
(G,\mathcal{O}_G)\stackrel{\delta_g\times \id}{\longrightarrow}
(G,\mathcal{O}_G)\times (G,\mathcal{O}_G)\stackrel{\mu}{\to}
(G,\mathcal{O}_G).
$$
The morphism $l_g$ is called the {\it left translation}. The {\it
right translation} $r_g$, $g\in G$, can be defined similarly. Denote
by $\overline{l}_g$ the following composition:
$$
(M,\mathcal{O}_M)=(g,\mathbb{C})\times
(M,\mathcal{O}_M)\stackrel{\delta_g\times \id}{\longrightarrow}
(G,\mathcal{O}_G)\times (M,\mathcal{O}_M)\stackrel{\nu}{\to}
(M,\mathcal{O}_M).
$$
The representations of $G_x$ in $T_x(M,\mathcal{O}_M)_{\bar 0}$ and
$T_x(M,\mathcal{O}_M)_{\bar 1}$  given by $G_x\ni h\mapsto (\d
\overline{l}_h)_x$ are called the {\it even and odd isotropy
representation}, respectively.

Assume that $(M,\mathcal{O}_M)$ is $(G,\mathcal{O}_G)$-homogeneous,
then $(\d \nu_x)_e$ is surjective. Hence,
$T_x(M,\mathcal{O}_M)\simeq \mathfrak{g}/\mathfrak{g}_x$. Denote by
$\Ad_G$ the adjoint representation of $G$ on $\mathfrak{g}$. Recall
that this representation is defined by
 $\Ad_G(g)(X)=(\d l_g\circ r_g)_e(X)$.
Clearly, $\Ad_G(h)$ transforms $\mathfrak{g}_x$ into itself for all
$h\in H$. It follows that there is a representation
$\widehat{\Ad}_G$ of $H$ in $\mathfrak{g}/\mathfrak{g}_x$ given by
$$
\widehat{\Ad}_G(h)(X+\mathfrak{g}_x)=
\Ad_G(h)(X)+\mathfrak{g}_x,\,\,X\in \mathfrak{g},\,\, h\in H.
$$
As in the classical case we have.

\medskip

\lem\label{isotropy} {\it The representation $\widehat{\Ad}_G$ is
equivalent to the isotropy representation and $(\d \nu_x)_e$
determines the corresponding equivalence. More precisely,
$\widehat{\Ad}_G|_{(\mathfrak{g}_x)_{\bar 0}}$ is equivalent to the
even isotropy representation and
$\widehat{\Ad}_G|_{(\mathfrak{g}_x)_{\bar 1}}$ to the odd one.}

\medskip
\noindent{\it Proof.} It is sufficient to check that for every $h\in
H$, the following diagram is commutative:
$$
\begin{CD}
T_e(G,\mathcal{O}_G)& @>{(\d \nu_x)_e}>> & T_x(M,\mathcal{O}_M)\\
@V{\Ad_G(h)}VV && @VV{(\d \overline{l}_h)_e}V\\
T_e(G,\mathcal{O}_G)& @>{(\d \nu_x)_e}>> & T_x(M,\mathcal{O}_M)
\end{CD}.
$$
It is easy to see that
$$
\nu_x\circ r_h=\nu_{hx}=\nu_x,\,\,\,\, \nu_x\circ
l_h=\overline{l}_h\circ \nu_x \quad \text{for all}\,\,\, h\in H.
$$
Therefore,
$$
\nu_x\circ r_h \circ l_h=\nu_x\circ l_h=\overline{l}_h\circ \nu_x.
\quad\Box
$$

\medskip

Let $(M,\mathcal{O}_M)$ be a $(G,\mathcal{O}_G)$-homogeneous
supermanifold, then $G$ acts on $\mathcal{T}_M$ by
\begin{equation}\label{action of G}
v\mapsto (\overline{l}_g^{-1})^* \circ v\circ (\overline{l}_g)^*.
\end{equation}

\medskip

\begin{center}
{\it 1.3 Split supermanifolds}
\end{center}

Let us describe the category $\verb"SSM"$ (split supermanifolds),
which was introduced in \cite{V}. Recall that a supermanifold
$(M,\mathcal{O}_M)$ is called {\it split} if $\mathcal{O}_M\simeq
\bigwedge_{\mathcal{F}_M}\mathcal{E}_M$ for a certain locally free
sheaf $\mathcal{E}_M$ over $M$.  We put
$$
\begin{array}{rl}
\operatorname{Ob}\, \verb"SSM"=\{ (M,\bigwedge
\mathcal{E}_M)\,\,|\,\, \mathcal{E}_M\,\,\text{ is a locally free
sheaf on}\,\, M\}.
\end{array}
$$
Equivalently, we can say that $\operatorname{Ob}\, \verb"SSM"$
consists of all split supermanifolds $(M,\mathcal{O}_M)$ with a
fixed isomorphism $\mathcal{O}_M\simeq \bigwedge \mathcal{E}_M$ for
a certain locally free sheaf $\mathcal{E}_M$ on $M$. Note that
$\mathcal{O}_M$ is naturally $\mathbb{Z}$-graded by
$(\mathcal{O}_M)_p\simeq \bigwedge^p \mathcal{E}_M$. All the time we
will consider this $\mathbb{Z}$-grading for elements from
$\operatorname{Ob}\, \verb"SSM"$. Further, if $X,Y\in
\operatorname{Ob}\, \verb"SSM"$, we put
$$
\begin{array}{rl}
\Hom(X,Y)=&\text{all morphisms of $X$ to $Y$}\\
 &\text{preserving the
$\mathbb{Z}$-gradings}.
\end{array}
$$

As in the category of supermanifolds, we can define in $\verb"SSM"$
a group object (split Lie supergroup), an action  of a split Lie
supergroup on a split supermanifold (split action) and a homogeneous
split supermanifold.

There is a functor $\gr$ from the category of supermanifolds to the
category of split supermanifolds. Let us briefly describe this
construction (see, e.g., \cite{Man,onipi}). Let $(M,\mathcal{O}_M)$
be a supermanifold. As above, denote by $\mathcal{J}_M\subset
\mathcal{O}_M$ the subsheaf of ideals generated by odd elements of
$\mathcal{O}_M$. Then by definition $\gr(M,\mathcal{O}_M)$ is the
split supermanifold $(M, \gr\mathcal{O}_M)$, where
$$
\gr\mathcal{O}_M= \bigoplus_{p\geq 0} (\gr\mathcal{O}_M)_p,\quad
\mathcal{J}_M^0:=\mathcal{O}_M, \quad (\gr\mathcal{O}_M)_p:=
\mathcal{J}_M^p/\mathcal{J}_M^{p+1}.
$$
In this case $(\gr\mathcal{O}_M)_1$ is a locally free sheaf and
there is a natural isomorphism of $\gr\mathcal{O}_M$ onto $\bigwedge
(\gr\mathcal{O}_M)_1$. If
$\psi=(\psi_{\red},\psi^*):(M,\mathcal{O}_M)\to (N,\mathcal{O}_N)$
is a morphism, then $\gr(\psi)=(\psi_{\red},\gr(\psi^*))$, where
$\gr(\psi^*):\gr \mathcal{O}_N \to \gr \mathcal{O}_M$ is defined by
$$
\gr(\psi^*)(f+\mathcal{J}_N^p): = \psi^*(f)+\mathcal{J}_M^p
\,\,\text{for}\,\, f\in (\mathcal{J}_N)^{p-1}.
$$
Recall that by definition every morphism $\psi$ of supermanifolds is
even and as a consequence sends $\mathcal{J}_N^p$ into
$\mathcal{J}_M^p$.

Let $(G,\mathcal{O}_G)$ be a Lie supergroup with the group morphisms
$\mu$, $\iota$ and $\varepsilon$. Then it is easy to see that $\gr
(G,\mathcal{O}_G)$ is a split Lie supergroup with the group
morphisms $\gr(\mu)$, $\gr(\iota)$ and $\gr(\varepsilon)$.
Similarly, an action $\nu:(G,\mathcal{O}_G)\times
(M,\mathcal{O}_M)\to (M,\mathcal{O}_M)$ gives rise to the action
$\gr(\nu):\gr(G,\mathcal{O}_G)\times \gr(M,\mathcal{O}_M)\to
\gr(M,\mathcal{O}_M)$.



Let $(M,\mathcal{O}_M)=(M,\bigwedge \mathcal{E}_M)$ be a split
supermanifold. Then the sheaf $\mathcal{T}_M$ is a
$\mathbb{Z}$-graded sheaf of Lie superalgebras. The
$\mathbb{Z}$-grading is given by
\begin{equation}\label{Z-grding of T}
(\mathcal{T}_M)_q=\{ v\in \mathcal{T}_M\,\mid\, v(\bigwedge^p
\mathcal{E}_M)\subset \bigwedge^{p+q} \mathcal{E}_M, \,\,p\geq 0\}.
\end{equation}
 The sheaf
$(\mathcal{T}_M)_q$, $q\in \mathbb{Z}$, is a locally free sheaf of
$\mathcal{F}_M$-modules. We will use the notation
$\mathfrak{v}(M,\mathcal{O}_M)_q:= H^0((\mathcal{T}_M)_q)$.

 It was shown in \cite{OniTransit} that
$\mathcal{E}_M^*\simeq (\mathcal{T}_M)_{-1}$. This isomorphism
identifies any sheaf homomorphism $\mathcal{E}_M\to \mathcal{F}_M$
with a derivation of degree $-1$ that is zero on $\mathcal{F}_M$.
Denote by $\mathbb{E}$ the vector bundle corresponding to
$\mathcal{E}_M$ and by $\mathbb{T}_{-1}$ the vector bundle
corresponding to $(\mathcal{T}_M)_{-1}$. It is easy to see that
$(\mathbb{T}_{-1})_x=T_x(M,\mathcal{O}_M)_{\bar 1}$, $x\in M$.

Assume in addition that $(M,\mathcal{O}_M)$ is
$(G,\mathcal{O}_G)$-homogeneous and the action of
$(G,\mathcal{O}_G)$ on $(M,\mathcal{O}_M)$ is split. Then the action
of $G$ on $\bigwedge \mathcal{E}_M$ given by $g\mapsto
\overline{l}_g$, $g\in G$, preserves the $\mathbb{Z}$-grading. Hence
the vector bundles $\mathbb{E}$ and $\mathbb{E}^*$ are
$G$-homogeneous. Furthermore, the corresponding action of $G$ on
$\mathcal{T}_M$ given by (\ref{action of G}) preserves the
$Z$-grading (\ref{Z-grding of T}).
 The following Lemma is
well-known.

\medskip

\lem\label{T_x isom E_x} {\it Assume that $(M,\mathcal{O}_M)$ is
$(G,\mathcal{O}_G)$-homogeneous and the action $\nu$ of
$(G,\mathcal{O}_G)$ on $(M,\mathcal{O}_M)$ is split. Then
$\mathbb{E}\simeq \mathbb{T}_{-1}^*$ as homogeneous vector bundles.
In particular,  $T_x(M,\mathcal{O}_M)_{\bar 1}\simeq \mathbb{E}_x^*$
as $G_x$-modules, where $G_x$ is the stabilizer of $x$ by the action
$\nu_{\red}$ of $G$ on $M$.

}

\medskip

\medskip

\begin{center}
   {\bf  2. Holomorphic functions on a complex
homogeneous supermanifold with compact reduction}
\end{center}

\begin{center}
{\it 2.1  The retract of a homogeneous supermanifold }
\end{center}

\medskip

Let $\mathfrak{g}_{\bar 0}$ be a Lie algebra and $V$ a
$\mathfrak{g}_{\bar 0}$-module. Denote by
$[\,,\,]_{\mathfrak{g}_{\bar 0}}$ the Lie bracket in
$\mathfrak{g}_{\bar 0}$ and by $\cdot$ the module operation in $V$.
We can construct a Lie superalgebra $\mathfrak{g}=\mathfrak{g}_{\bar
0}\oplus \mathfrak{g}_{\bar 1}$ putting $\mathfrak{g}_{\bar 1}=V$
and defining the Lie bracket by the following formula:
\begin{equation}\label{bracket split}
[X,Y]=\left\{
        \begin{array}{ll}
          [X,Y]_{\mathfrak{g}_{\bar 0}}, & \hbox{if $X,Y\in \mathfrak{g}_{\bar 0}$;} \\
          X\cdot Y, & \hbox{if $X\in \mathfrak{g}_{\bar 0}$ and $Y\in V$;} \\
          0, & \hbox{if $X,Y\in V$.}
        \end{array}
      \right.
\end{equation}

Let $G$ be a Lie group, $\mathfrak{g}_{\bar 0} = \Lie  G$, and $V$ a
$G$-module. Assume that $\mathfrak{g}=\mathfrak{g}_{\bar 0}\oplus V$
is the Lie superalgebra with the Lie bracket given by (\ref{bracket
split}). Using this data we can construct a Lie supergroup in the
following way. Let us describe its Harish-Chandra pair. Denote by
$\alpha_G$ the representation of $G$ on $\mathfrak{g}$ given by:
\begin{equation}\label{alpha_G, split Lie sgroup}
\begin{array}{c}
\alpha_G| \mathfrak{g}_{\bar 0} =\,\text{adjoint representation
of $G$ on $\mathfrak{g}_{\bar 0}$},\\
(\alpha_G| \mathfrak{g}_{\bar 1})(g)(v):=g\cdot v,\, g\in G,\,v\in
V\, (\text{the given module operation}).
\end{array}
\end{equation}
Now the Harish-Chandra pair $(G,\mathfrak{g})$ is well-defined.

Let $(G,\mathcal{O}_G)$ be a Lie supergroup, $\mathfrak{g}=\Lie
(G,\mathcal{O}_G)$. Denote by $\mathfrak{g}'$ the Lie superalgebra
such that $\mathfrak{g}'\simeq \mathfrak{g}$ as vector superspaces,
the Lie bracket is defined by (\ref{bracket split}) with
$V:=\mathfrak{g}_{\bar 1}$ and $X\cdot Y=[X,Y]$ for $X\in
\mathfrak{g}_{\bar 0}$, $Y\in \mathfrak{g}_{\bar 1}$.

\medskip

\t\label{gr(G po H)isom grG po gr H} {\it $1.$ The Lie supergroup
$\gr(G,\mathcal{O}_G)$ is determined by the Harish-Chandra pair
$(G,\mathfrak{g}')$.

$2.$ If $(H,\mathcal{O}_H)$ is a closed subsupergroup of
$(G,\mathcal{O}_G)$, then
$$
\gr((G,\mathcal{O}_G)/(H,\mathcal{O}_H))\simeq
\gr(G,\mathcal{O}_G)/\gr (H,\mathcal{O}_H).
$$

$3.$ The supermanifold $\gr(G,\mathcal{O}_G)/ \gr(H,\mathcal{O}_H)$
is split. If $\mathbb{E}$ is the corresponding homogeneous bundle,
then it is determined by the $H$-module $(\mathfrak{g}_{\bar
1}/\mathfrak{h}_{\bar 1})^*$. }

\medskip

\noindent{\it Proof.} To prove the first statement of the theorem,
we have to prove that
\begin{equation}\label{komutator split}
[X,Y]_{\mathfrak{g}'}=\left\{
                        \begin{array}{ll}
                          [X,Y]_{\mathfrak{g}}, & \hbox{if $X,Y\in \mathfrak{g}'_{\bar 0}$ or $X\in \mathfrak{g}'_{\bar 0}$,
$Y\in \mathfrak{g}'_{\bar 1}$;} \\
                          0, & \hbox{if $X,Y\in \mathfrak{g}'_{\bar 1}$.}
                        \end{array}
                      \right.
\end{equation}
 Here $[\,\,,\,]_{\mathfrak{g}}$ and $[\,\,,\,]_{\mathfrak{g}'}$ are the
Lie brackets in $\mathfrak{g}$ and $\mathfrak{g}'$, respectively.
Let us take $X_e,Y_e\in
T_e(\gr(G,\mathcal{O}_G))=T_e(G,\mathcal{O}_G)$. We put
$$
\begin{array}{ll}
X=(X_e\otimes \id)\circ \mu^*, &X'=(X_e\otimes \id)\circ (\gr
\mu)^*,\\
Y=(Y_e\otimes \id)\circ \mu^*, &Y'=(Y_e\otimes \id)\circ (\gr
\mu)^*,\\
Z=[X,Y]_{\mathfrak{g}}, &Z'=[X',Y']_{\mathfrak{g}'}.
\end{array}
$$
To prove (\ref{komutator split}) it is enough to show that
$\delta_e\circ Z= \delta_e\circ Z'$ if $X,Y\in \mathfrak{g}'_{\bar
0}$ or $X\in \mathfrak{g}'_{\bar 0}$, $Y\in \mathfrak{g}'_{\bar 1}$,
and that $\delta_e\circ Z'=0$ if  $X,Y\in \mathfrak{g}'_{\bar 1}$.

Let us take $f\in (\mathcal{O}_G)_p$, then
$\mu^*(f)=g_p+g_{p+2}+\ldots$ and $(\gr \mu)^*(f)=g_p$, where
$g_i\in (\mathcal{O}_{G\times G})_i$. It is easy to see that $g_p\ne
0$ (it follows, for example, from the identity axiom of a Lie
supergroup). Further, using (\ref{left inv vect field}) we get
$$
\begin{array}{c}
 \delta_e\circ Z= \delta_e\circ((-1)^{p(X)p(Y)}(
Y_e\otimes X_e \otimes \id)- ( X_e\otimes
Y_e \otimes \id))\circ ((\id \times \mu)\circ\mu)^*\\
=((-1)^{p(X)p(Y)}(Y_e\otimes X_e)-(X_e\otimes Y_e))\circ \mu^*.
\end{array}
$$
Similarly,
$$
\begin{array}{c}
 \delta_e\circ Z'=
((-1)^{p(X)p(Y)}(Y_e\otimes X_e)-(X_e\otimes Y_e))\circ (\gr\mu)^*.
\end{array}
$$
Assume that $X_e,Y_e\in T_e(G,\mathcal{O}_G)_{\bar 0}$ and $f\in
(\mathcal{O}_G)_p$. Then
$$
\begin{array}{c}
\delta_e\circ Z(f)= ((Y_e\otimes X_e)-(X_e\otimes
Y_e))\circ \mu^*(f)=\\
((Y_e\otimes X_e)-(X_e\otimes
Y_e))(g_p+g_{p+2}+\ldots)=\\
\left\{
  \begin{array}{ll}
    ((Y_e\otimes X_e)- (X_e\otimes
Y_e))(g_0), & \hbox{if $p=0$;} \\
    0, & \hbox{if $p\geq 1$.}
  \end{array}
\right.
\end{array}
$$
Similarly,
$$
\begin{array}{c}
\delta_e\circ Z'(f)= ((Y_e\otimes X_e)-(X_e\otimes
Y_e))\circ (\gr\mu)^*(f)=\\
((Y_e\otimes X_e)-(X_e\otimes Y_e))(g_p)= \left\{
  \begin{array}{ll}
    ((Y_e\otimes X_e)- (X_e\otimes
Y_e))(g_0), & \hbox{if $p=0$;} \\
    0, & \hbox{if $p\geq 1$.}
  \end{array}
\right.
\end{array}
$$
Hence, in this case $\delta_e\circ Z=\delta_e\circ Z'$.

Assume that $X_e\in T_e(G,\mathcal{O}_G)_{\bar 0}$, $Y_e\in
T_e(G,\mathcal{O}_G)_{\bar 1}$ and $f\in (\mathcal{O}_G)_p$. Then as
above we get
$$
\begin{array}{c}
\delta_e\circ Z(f)= \left\{
  \begin{array}{ll}
    ((Y_e\otimes X_e)- (X_e\otimes
Y_e))(g_1), & \hbox{if $p=1$;} \\
    0, & \hbox{if $p=0$ or $p\geq 2$.}
  \end{array}
\right.
\end{array}
$$
and
$$
\begin{array}{c}
\delta_e\circ Z'(f)= \left\{
  \begin{array}{ll}
    ((Y_e\otimes X_e)- (X_e\otimes
Y_e))(g_1), & \hbox{if $p=1$;} \\
    0, & \hbox{if $p=0$ or $p\geq 2$,}
  \end{array}
\right.
\end{array}
$$
Hence, in this case $\delta_e\circ Z=\delta_e\circ Z'$ as well.

Assume that $X_e,Y_e\in T_e(G,\mathcal{O}_G)_{\bar 1}$. Then
$$
\begin{array}{c}
\delta_e\circ Z(f)= ((Y_e\otimes X_e)+(X_e\otimes Y_e))\circ
(\gr\mu)^*(f)=\\
 ((Y_e\otimes X_e)-(X_e\otimes Y_e))(g_p)=0, \,\,\,p\geq 0.
\end{array}
$$
The proof of (\ref{komutator split}) is complete.

To prove the second statement of the theorem, denote by $\nu$ the
action of the Lie supergroup $(G,\mathcal{O}_G)$ on the
supermanifold $(M,\mathcal{O}_M):=
(G,\mathcal{O}_G)/(H,\mathcal{O}_H)$. It is easy to see that the
action $\gr \nu$ is transitive on $\gr(M,\mathcal{O}_M)$ (see
\cite[Lemma $5$]{V}). Hence, it is enough to show that the
stabilizer of the point $eH\in G/H$ is $\gr (H,\mathcal{O}_H)$. Note
that $H$ is the stabilizer of $eH$ by the action $\nu_{\red}=(\gr
\nu)_{\red}$ and $\Ker (\d\gr\nu_{eH})_e=\Ker (\d\nu_{eH})_e$ as
vector spaces. Now this assertion follows from the first one.

To complete the proof, note that $\gr \gr(G,\mathcal{O}_G)=
\gr(G,\mathcal{O}_G)$. The last assertion follows from Lemmas
\ref{isotropy} and \ref{T_x isom E_x}.$\Box$

\medskip

The third part of this theorem was proved in a different way in
\cite[Theorem $2$]{V}.

The description of the vector bundle determined by a split
homogeneous supermanifold gives the following proposition.

\medskip
\l\label{bundel of split homogen smf} {\it Assume that $M=G/H$ is a
complex compact homogeneous manifold and $(G,\mathcal{O}_G)$ is a
Lie supergroup. Denote by $\mathbb{E}$ the vector bundle determined
by a split $(G,\mathcal{O}_G)$-homogeneous supermanifold
$(M,\bigwedge \mathcal{E})$. Then $\mathbb{E}$ is a homogeneous
subbundle of a trivial homogeneous bundle $\mathbb{V}$ determined by
a certain $G$-module $V$.

Conversely, any homogeneous subbundle $\mathbb{E}$ of a trivial
homogeneous bundle $\mathbb{V}$ determined by a $G$-module $V$
corresponds to a certain split homogeneous supermanifold
$(M,\bigwedge \mathcal{E})$.

}

\medskip

\noindent{\it Proof.} Assume that $(M,\mathcal{O}_M)$ is
$(G,\mathcal{O}_G)$-homogeneous. Then $G$ acts on
$(M,\mathcal{O}_M)$ by left translations $\overline{l}_g$, $g\in G$.
Note that $\overline{l}_g$ is not an automorphism of the vector
bundle $\mathbb{E}$. By Theorem \ref{gr(G po H)isom grG po gr H} we
have $(M,\mathcal{O}_M)\simeq \gr (G,\mathcal{O}_G)/\gr
(H,\mathcal{O}_H)$. The left translations determined by the action
of $\gr (G,\mathcal{O}_G)$ are automorphisms of $\mathbb{E}$. Hence,
$\mathbb{E}$ is a $G$-homogeneous vector bundle. (The fact that the
vector bundle determined by a split homogeneous supermanifold is
homogeneous, was also noticed in \cite{OniTransit}.) Furthermore,
the vector space $V: = \mathfrak{v}(M,\mathcal{O}_M)_{-1}^*$ is a
finite dimensional $G$-module because $M$ is compact. Since
$(M,\mathcal{O}_M)$ is homogeneous, the $H$-equivariant map
$$
\ev_{eH}:\mathfrak{v}(M,\mathcal{O}_M)_{-1}\to
T_{eH}(M,\mathcal{O}_M)\simeq \mathbb{E}_{eH}^*
$$
 is surjective.
Hence, the dual map $\mathbb{E}_{eH}\to V$ is injective and
$\mathbb{E}$ is a homogeneous subbundle of the trivial bundle
$\mathbb{V}$ determined by the $H$-module $V$.

Conversely, we put $E= \mathbb{E}_{eH}$. This is an $H$-module.
Denote by $(G,\mathcal{O}_G)$ the Lie supergroup determined by the
Harish-Chandra pair $(G,\mathfrak{g})$, where
$\mathfrak{g}=\mathfrak{g}_{\bar 0}\oplus V^*$, $\mathfrak{g}_{\bar
0}=\Lie  G$ and the Lie bracket is given by (\ref{bracket split}).
Let also $(H,\mathcal{O}_H)$ be the Lie subsupergroup of
$(G,\mathcal{O}_G)$ determined by the Harish-Chandra subpair
$(H,\mathfrak{h})$, where $\mathfrak{h}=\mathfrak{h}_{\bar 0}\oplus
E'$, $E'=\Ker(V^*\to E^*)$ and $\mathfrak{h}_{\bar 0}=\Lie H$. Then
by Theorem \ref{gr(G po H)isom grG po gr H} the homogeneous
supermanifold $(G/H,\mathcal{O}_{G/H})$ is split and the
corresponding homogeneous vector bundle $\mathbb{E}$ is determined
by the $H$-module $(V^*/E')^*\simeq E$. The proof is complete.$\Box$



\medskip

\lem\label{H^0(E)=0} {\it Let $G$ be a complex Lie group, $H\subset
G$ a closed complex Lie subgroup, $V$ a $G$-module and $E\subset V$
an $H$-submodule. Assume that $G/H$ is compact and connected. Denote
by $\mathbb{E}$ the homogeneous vector bundle that corresponds to
$E$. The following conditions are equivalent:
\begin{enumerate}
  \item[1)] non-trivial $G$-modules $W$ such that $W\subset
E$ do not exist;
  \item[2)] $\Gamma(\mathbb{E})=\{0\}$.
\end{enumerate}
}

\medskip

\noindent{\it Proof.} $1)\Rightarrow 2)$ Assume that
$\Gamma(\mathbb{E})\ne \{0\}$. Denote by $\mathbb{V}$ the
homogeneous vector bundle determined by the $H$-module $V$. It is
trivial and the evaluation map $\Gamma(\mathbb{V})\to V$, $s\mapsto
s_x$, $x=eH\in G/H$, is an isomorphism of $G$-modules. The bundle
$\mathbb{E}$ is a subbundle of $\mathbb{V}$, hence there is an
inclusion $\gamma:\Gamma(\mathbb{E}) \to \Gamma(\mathbb{V})$. Denote
by $W$ the image of $\gamma(\Gamma(\mathbb{E}))$ in $V$ by the map
$s\mapsto s_x$. It is a $G$-submodule in $V$ and it is non-trivial
by the assumption. Consider the commutative diagram:
$$
\begin{CD}
\Gamma(\mathbb{E}) &@>>> & \Gamma(\mathbb{V})\\
@VVV & & @VVV\\
E & @>>> & V
\end{CD}
$$
where the horizontal arrows are inclusions and the vertical arrows
are evaluation maps at the point $x$. We see that the image of $E$
contains $W$. We arrive at a contradiction.

$2)\Rightarrow 1)$ Assume that there is a non-trivial $G$-module $W$
in $E$. Then there is a trivial subbundle $\mathbb{W}$ in
$\mathbb{E}$, where $\mathbb{W}$ is the homogeneous vector bundle,
determined by the $H$-module $W$. Hence $\Gamma(\mathbb{E})\ne
\{0\}$.$\Box$

\medskip

\begin{center}
   {\it 2.2 Odd fundamental vector fields on a split homogeneous supermanifold.}
\end{center}

Let $(G,\mathcal{O}_G)$ be a complex Lie supergroup. It was proved
in \cite{V} that it is isomorphic to the Lie supergroup
$(G,\widehat{\mathcal{O}}_G)$ determined by the Harish-Chandra pair
$(G,\mathfrak{g})$ using Koszul construction, see $1.3$. The
isomorphism is given by the following formula:
\begin{equation}\label{isomorphism of supergroups}
\Phi : \mathcal{O}_G \to\widehat{\mathcal{O}}_G,\,\,\,
 \Phi(f)(X)(g)= (-1)^{p(X)p(f)} (X(f))_{\red}(g).
\end{equation}
We will identify the Lie supergroups $(G,\mathcal{O}_G)$ and
$(G,\widehat{\mathcal{O}}_G)$ using this isomorphism. The sheaf
$\mathcal{O}_G$ is $\mathbb{Z}$-graded, this $\mathbb{Z}$-grading is
induced by the following $\mathbb{Z}$-grading:
\begin{equation}\label{grading}
\Hom(\bigwedge \mathfrak{g}_{\bar 1}, \mathcal{F}_G)=
\bigoplus_{q\geq 0} \Hom(\bigwedge^q \mathfrak{g}_{\bar 1},
\mathcal{F}_G).
\end{equation}
In other words, the Lie supergroup $(G,\mathcal{O}_G)$ possesses a
global odd coordinate system. Namely, let $(\xi_i)$ be a basis of
$\mathfrak{g}_{\bar 1}$. Let $f^{\xi_i}\in \mathcal{O}_G$ so
$f^{\xi_i}\circ \gamma \in \Hom(\mathfrak{g}_{\bar 1},
\mathcal{F}_G)$ and $f^{\xi_i}\circ \gamma(\xi_j)=\delta_{ij}$. Then
$(f^{\xi_i})$ is a global odd coordinate system on
$(G,\mathcal{O}_G)$.

Our aim is now to describe right invariant vector fields in the
chosen odd coordinates for a split Lie supergroup. Let us take $X\in
\mathfrak{g}$, $Y\in \mathfrak{U}(\mathfrak{g})$ and $g\in G$. Using
(\ref{isomorphism of supergroups}) we get
$$
(X(f))(Y)(g)=(-1)^{p(Y)p(X(f))} (Y\cdot X(f))_{\red}(g)=(-1)^{p(X)}
f(Y\cdot X)(g).
$$
Hence, the right invariant vector field $X$ is determined by
\begin{equation}\label{right invar vect field}
(X(f))(Y)(g)= (-1)^{p(X)} f(Y\cdot X)(g).
\end{equation}

Assume in addition that $(G,\mathcal{O}_G)$ is a split Lie
supergroup. From Theorem~\ref{gr(G po H)isom grG po gr H} it follows
that this is equivalent to $[\mathfrak{g}_{\bar
1},\mathfrak{g}_{\bar 1}]=\{0\}$. In this case, the map $\gamma$
from (\ref{isomorphism}) is a homomorphism of algebras. Denote by
$X_{\xi_i}$ the right invariant vector field which corresponds to
$\xi_i$. By (\ref{right invar vect field}) we get for $Y\in
\bigwedge^p \mathfrak{g}_{\bar 1}$ and $g\in G$,
\begin{equation}\label{X_xi_i}
\begin{split}
(X_{\xi_i}(f^{\xi_j})\circ \gamma)(Y)(g)= (X_{\xi_i}(f^{\xi_j}))
(Y)(g)= - f^{\xi_j}(Y\cdot \xi_i)=\\
\left\{
                                    \begin{array}{ll}
                                      0, & \hbox{$p\ne 0$;} \\
                                      -\delta_{ij}, & \hbox{$p=0$.}
                                    \end{array}
                                  \right.
\end{split}
\end{equation}
Let $\mu$ be the multiplication morphism of $(G,\mathcal{O}_G)$.
Since $\gr \mu=\mu$, we get that $X_{\xi_i}\in
\mathfrak{v}(G,\mathcal{O}_G)_{-1}$ by (\ref{left inv vect field}).
It follows that $X_{\xi_i}$ is completely determined by
(\ref{X_xi_i}) and has the form $-\frac{\partial}{\partial
f^{\xi_i}}$ in the chosen odd coordinates. We have proved the
following result:

\medskip

\lem\label{vector fields} {\it Let $(G,\mathcal{O}_G)$ be a split
Lie supergroup and $(f^{\xi_i})$ the global odd coordinate system
described above. Then the vector fields $\frac{\partial}{\partial
f^{\xi_i}}$, $i=1,\ldots, \dim \mathfrak{g}_{\bar 1}$ are right
invariant.

 }

\medskip

Now we are able to prove the following lemma.

\medskip

\lem\label{H^0(F)=0 follows H^0(wedge F)=0} {\it Let
$(M,\mathcal{O}_M)$ be a split homogeneous supermanifold and
$\mathcal{O}_M\simeq \bigwedge \mathcal{E}$. If
$H^0(\mathcal{E})=0$, then $H^0(\bigwedge^p\mathcal{E})=0$ for all
$p>0$. }

\medskip

\noindent{\it Proof.} By Theorem \ref{gr(G po H)isom grG po gr H} we
may assume that $(M,\mathcal{O}_M)$ is a
$(G,\mathcal{O}_G)$-homogeneous supermanifold, where
$(G,\mathcal{O}_G)$ is a split Lie supergroup. Let $\mathfrak{g}$ be
a Lie superalgebra of $(G,\mathcal{O}_G)$. Denote by
$(H,\mathcal{O}_H)$ the stabilizer of a point $x\in M$. Then
$(M,\mathcal{O}_M)\simeq (G/H,\mathcal{O}_{G/H})$ and
$\mathcal{O}_{G/H}\subset \mathcal{O}_{G}$, see $1.1$. In
\cite[Proposition $5$]{V} it was shown that if $[\mathfrak{g}_{\bar
1}, \mathfrak{g}_{\bar 1}]=\{0\}$, then there is an isomorphism of
sheaves $\bigwedge \mathcal{E}\to \mathcal{O}_{G/H}$ such that the
composition $\bigwedge \mathcal{E}\to
\mathcal{O}_{G/H}\hookrightarrow \mathcal{O}_{G}$ preserves the
$\mathbb{Z}$-gradings of sheaves.

Assume that $H^0(\bigwedge^p\mathcal{E})\ne 0$ for a certain $p>0$
and let
 $f\in H^0(\bigwedge^p\mathcal{E})$, $f\ne 0$. Then $f\in
H^0((\mathcal{O}_G)_p) \simeq H^0(\mathcal{F}_G)\otimes \bigwedge^p
\mathfrak{g}_{\bar 1}^*$. Let $(f^{\xi_i})$ be the global odd
coordinate system described above. Then we can write $f$ in the
following form:
$$
f=\sum_{i_1<\cdots< i_p}
f_{i_1,\ldots,i_p}f^{\xi_{i_1}}\wedge\cdots\wedge f^{\xi_{i_p}},
\,\,\,f_{i_1,\ldots,i_p}\in H^0(\mathcal{F}_G).
$$
Without loss of generality, we may assume that $f^{\xi_1}$ occurs on
the right-hand side of this equation. Hence, we can write
$f=f^{\xi_1} g + h$, where $0\ne g\in H^0(\mathcal{F}_G)\otimes
\bigwedge^{p-1} \mathfrak{g}_{\bar 1}^*$ and $h$ does not depend on
$f^{\xi_1}$. By Lemma \ref{vector fields}, the vector field
$\frac{\partial}{\partial f^{\xi_1}}$ is right invariant. It follows
that
$$
\frac{\partial}{\partial f^{\xi_1}}(f)=g\in
(\mathcal{O}_{G/H})_{p-1} \simeq H^0(\bigwedge^{p-1}\mathcal{E}).
$$
Hence, $H^0(\bigwedge^{p-1}\mathcal{E})\ne 0$. By induction, the
proof is complete.$\Box$

\medskip

We need the following lemma.

\medskip

\lem\label{funk_retr} {\it Let $(M,\mathcal O_M)$ be a
supermanifold. If $H^0(\gr\mathcal O_M)\simeq\Bbb C$ then
$H^0(\mathcal O_M)\simeq\Bbb C$.}

\medskip

\noindent{\it Proof.} If $H^0(\gr\mathcal O_M )\simeq\Bbb C$,
then $H^0((\gr\mathcal O_M)_p)=\{0\}$ for all $p > 0$ and 
$H^0((\gr\mathcal O_M)_0 )\simeq\Bbb C$. For each $p\ge 0$, we have
the exact sequence

$$
0\to H^0(\mathcal J_M^{p+1})\to H^0(\mathcal J_M^p)\to
H^0((\gr\mathcal O_M)_p),
$$
where $\mathcal J_M$ is the sheaf generated by odd elements of
 $\mathcal O_M$.
We have $H^0(\mathcal J_M^p)=\{0\}$ if $p$ is sufficiently large.
Using induction, we see that $H^0(\mathcal J_M^p)=\{0\}$ for all
$p>0$. For $p=0$, we have the exact sequence
$$
0\to H^0(\mathcal J_M^{0})= H^0(\mathcal O_M)\to H^0((\gr\mathcal
O_M)_0)\simeq \mathbb{C}.
$$
Obviously, $H^0(\mathcal O_M)\supset \Bbb C$, since there are
constant functions on every supermanifold. Hence, $H^0(M,\mathcal
O_M)\simeq\Bbb C$.$\Box$

\medskip

\begin{center}
{\it 2.3 The main result}
\end{center}

\medskip

\t\label{funk1} {\it Let $(M,\mathcal O_M)$ be a $(G,\mathcal
O_G)$-homogeneous supermanifold, $M$ a compact connected manifold,
$(H,\mathcal O_H)$ the stabilizer of a point $x\in M$,
$\mathfrak{g}=\Lie (G,\mathcal O_G)$, $\mathfrak{h}=\Lie (H,\mathcal
O_H)$. Consider the exact sequence of $H$-modules:
$$
0\to \mathfrak{h}_{\bar 1}\to \mathfrak{g}_{\bar
1}\stackrel{\gamma}{\to} \mathfrak{g}_{\bar 1}/\mathfrak{h}_{\bar
1}\to 0.
$$
If there do not exist non-trivial $G$-modules $W\subset
\mathfrak{g}_{\bar 1}^*$ such that $W\subset \Im \gamma^*$, then
$H^0(\mathcal{O}_M)\simeq \mathbb{C}$. If in addition $(M,\mathcal
O_M)$ is split, then the converse statement is also true.

}

\medskip

\noindent{\it Proof.} Assume that $(M,\mathcal O_M)$ is split. By
Theorem \ref{gr(G po H)isom grG po gr H} we have $(M,\mathcal
O_M)\simeq \gr (G,\mathcal O_G)/\gr (H,\mathcal O_H)$. Denote by
$\mathbb{E}$ the vector bundle determined by $(M,\mathcal O_M)$.
 Put $\mathfrak{g}'=\Lie  (gr
(G,\mathcal{O}_G))$ and $\mathfrak{h}'=\Lie  (gr
(H,\mathcal{O}_H))$. By Lemmas \ref{isotropy} and \ref{T_x isom
E_x}, we see that $\mathbb{E}_x\simeq (\mathfrak{g}'_{\bar
1}/\mathfrak{h}'_{\bar 1})^*$ as $H$-modules. By Theorem \ref{gr(G
po H)isom grG po gr H}, we have that $\mathfrak{g}'_{\bar 1}=
\mathfrak{g}_{\bar 1}$ as $G$-modules and $\mathfrak{h}'_{\bar 1}=
\mathfrak{h}_{\bar 1}$ as $H$-modules, hence by assumption
non-trivial $G$-modules $W'\subset (\mathfrak{g}'_{\bar 1})^*$ such
that $W'\subset \Im \gamma'^*$, where
$\gamma'^*:(\mathfrak{g}'_{\bar 1}/\mathfrak{h}'_{\bar 1})^* \to
(\mathfrak{g}'_{\bar 1})^*$, do not exist. It follows from Lemma
\ref{H^0(E)=0} that $\Gamma(\mathbb{E})=\{0\}$. Furthermore, by
Lemma \ref{H^0(F)=0 follows H^0(wedge F)=0}, we get
$\Gamma(\bigwedge^p \mathbb{E})=\{0\}$ for all $p>0$. Hence,
$H^0(\mathcal O_M)=H^0(\mathcal{F}_M)\simeq \mathbb{C}$.

Conversely, if $H^0(\mathcal O_M) \simeq \mathbb{C}$ then
$H^0(\mathcal F_M) \simeq \mathbb{C}$ and $\Gamma(\bigwedge^p
\mathbb{E}) = \{0 \} $ for $p>0$. It follows from Lemma
\ref{H^0(E)=0} that non-trivial $G$-modules $W\subset
\mathfrak{g}_{\bar 1}^*$ such that $W\subset \Im \gamma^*$ do not
exist.

For non-split supermanifolds the assertion follows from Lemma
\ref{funk_retr} and Theorem \ref{gr(G po H)isom grG po gr H}.$\Box$

\medskip

\noindent{\bf Corollary.} {\it Let $(M,\mathcal O_M)$ be a
$(G,\mathcal O_G)$-homogeneous supermanifold, $M$ a compact
connected manifold, $(H,\mathcal O_H)$ the stabilizer of a point
$x\in M$, $\mathfrak{g}=\Lie (G,\mathcal O_G)$, $\mathfrak{h}=\Lie
(H,\mathcal O_H)$ and $\mathfrak{g}_{\bar 1}$ an irreducible
$G$-module. If the odd dimension of $(H,\mathcal{O}_H)$ is equal to
$0$, then $H^0(\mathcal O_M)\simeq \bigwedge (\mathfrak{g}_{\bar
1}^*)$. Otherwise, $H^0(\mathcal O_M)\simeq \mathbb{C}$. }

\medskip

The following proposition can be useful for practical applications:

\medskip

\l\label{funk_dop} {\it Let $(M,\mathcal O_M)$ be a $(G,\mathcal
O_G)$-homogeneous supermanifold, $M$ a compact connected manifold,
$(H,\mathcal O_H)$ the stationary subsupergroup of a point $x\in M$,
$\mathfrak{g}=\Lie  (G,\mathcal O_G)$, $\mathfrak{h}=\Lie
(H,\mathcal O_H)$. Assume that $\mathfrak{g}_{\bar 1}$ is a
completely reducible $G$-module. Consider the exact sequence of
$H$-modules:
$$
0\to \mathfrak{h}_{\bar 1}\stackrel{\delta}{\to}  \mathfrak{g}_{\bar
1}\stackrel{\gamma}{\to} \mathfrak{g}_{\bar 1}/\mathfrak{h}_{\bar
1}\to 0.
$$
Let $W\subset \Im \gamma^*$ be the maximal $G$-module and let
$Y=\{y\in \mathfrak{g}_{\bar 1}\,\mid\, W(y)=0 \}$. If
$\delta(\mathfrak{h}_{\bar 1}) \subset Y$, then
$H^0(\mathcal{O}_M)\simeq \bigwedge W$.
 If in addition $(M,\mathcal
O_M)$ is split, then $(M,\mathcal O_M)\simeq (N,\mathcal O_N)\times
(\pt,\bigwedge W)$, where $(N,\mathcal O_N)$ is a homogeneous
supermanifold such that $H^0(\mathcal{O}_N)\simeq \mathbb{C}$.

}

\medskip

Let us first prove the following lemma.

\medskip

\lem\label{product} {\it Let $(G_i,\mathcal{O}_{G_i})$, $i=1,2$, be
two Lie supergroups and $(H_i,\mathcal{O}_{H_i})\subset
(G_i,\mathcal{O}_{G_i})$  closed Lie subgroups. Then
$$
(G_1\times G_2/H_1\times H_2,\mathcal{O}_{G_1\times G_2/H_1\times
H_2})\simeq (G_1/H_1,\mathcal{O}_{G_1/H_1})\times
(G_2/H_2,\mathcal{O}_{G_2/H_2}).
$$

}

\medskip

\noindent{\it Proof.} Denote by $\nu_i$ the actions of
$(G_i,\mathcal{O}_{G_i})$ on $(G_i/H_i,\mathcal{O}_{G_i/H_i})$. Then
$\nu_1\times \nu_2$ is the action of $(G_1,\mathcal{O}_{G_1})\times
(G_2,\mathcal{O}_{G_2})$ on $(G_1/H_1,\mathcal{O}_{G_1/H_1})\times
(G_2/H_2,\mathcal{O}_{G_2/H_2})$. Let us compute the stabilizer of
the point $x=eH_1\times eH_2$. Denote by $(H',\mathfrak{h}')$ the
Harish-Chandra subpair determined by this stabilizer. Then $H'$ is
the stabilizer of $x$ corresponding to the action $(\nu_1\times
\nu_2)_{\red}=(\nu_1)_{\red}\times (\nu_2)_{\red}$. Hence $H' =
H_1\times H_2$. Furthermore,
$$
\mathfrak{h}'=\Ker( \d(\nu_1\times \nu_2)_x )_{e}= \Ker(
\d(\nu_1)_{eH_1} )_{e} \oplus \Ker( \d(\nu_2)_{eH_2}
)_{e}=\mathfrak{h}_1 \oplus \mathfrak{h}_2,
$$
where $\mathfrak{h}_i=\Lie  (H_i,\mathcal{O}_{H_i})$. It follows
that $(H',\mathfrak{h}')$ is the Harish-Chandra subpair of the Lie
group $(H_1,\mathcal{O}_{H_1}) \times (H_2,\mathcal{O}_{H_2})$. The
lemma follows.$\Box$

\medskip

\noindent{\it Proof of Proposition \ref{funk_dop}.} If $(M,\mathcal
O_M)$ is split, then by Theorem \ref{gr(G po H)isom grG po gr H} we
may assume that $[\mathfrak{g}_{\bar 1}, \mathfrak{g}_{\bar
1}]=\{0\}$. Let $V\subset \mathfrak{g}_{\bar 1}^*$ be a
$G$-submodule such that $\mathfrak{g}_{\bar 1}^*= W\oplus V$. Put
$X=\{x\in \mathfrak{g}_{\bar 1}\,\mid\, V(x)=0 \}$. The
subsuperspaces $\mathfrak{g}_{\bar 0}\oplus Y$ and $X$ are
subsupralgebras of $\mathfrak{g}$. Denote by
$(G_1,\mathcal{O}_{G_1})$ and by $(G_2,\mathcal{O}_{G_2})$ the Lie
subsupergroups of $(G,\mathcal{O}_{G})$ determined by the
Harish-Chandra subpairs $(G,\mathfrak{g}_{\bar 0}\oplus Y)$ and
$(e,X)$, respectively. Then $(G,\mathcal{O}_G)\simeq
(G_1,\mathcal{O}_{G_1})\times (G_2,\mathcal{O}_{G_2})$ and
$(H,\mathcal{O}_H)\subset (G_1,\mathcal{O}_{G_1})$. By Lemma
\ref{product}, we have
$$
(M,\mathcal O_M)\simeq (G_1/H,\mathcal{O}_{G_1/H})\times
(G_2,\mathcal{O}_{G_2})= (G_1/H,\mathcal{O}_{G_1/H})\times (\pt,
\bigwedge X^*).
$$
Hence, it is enough to show that $H^0(\mathcal{O}_{G_1/H}) \simeq
\mathbb{C}$. Let
$$
\begin{array}{l}
a:Y^*\to \mathfrak{g}_{\bar 1}^*, \,\,\, f\mapsto f\circ \pr_Y,\\
b:(Y/\mathfrak{h}_{\bar 1})^*\to (\mathfrak{g}_{\bar
1}/\mathfrak{h}_{\bar 1})^*, \,\,\, f\mapsto f\circ \pr_{\gamma(Y)}.
\end{array}
$$
Then $a$ is a homomorphism of $G$-modules, $b$ is a homomorphism of
$H$-modules, $a\circ (\gamma|_Y)^* = \gamma^* \circ b$ and $\Im
a=V$. If there is a non-trivial $G$-module in $\Im (\gamma|_Y)^*$,
then there is a non-trivial $G$-module in $\Im \gamma^*  \cap V$.
This contradicts the assumption that $W$ is maximal. By Theorem
\ref{funk1}, we get that $H^0(\mathcal{O}_{G_1/H}) \simeq
\mathbb{C}$.

Assume now that $(M,\mathcal O_M)$ is not split. The subsuperspace
$\mathfrak{g}_{\bar 0}\oplus Y\subset \mathfrak{g}$ is again a Lie
subsuperalgebra. Denote by $(G_1,\mathcal{O}_{G_1})$ the Lie
subsupergroup of $(G,\mathcal{O}_{G})$ determined by the
Harish-Chandra subpair $(G,\mathfrak{g}_{\bar 0}\oplus Y)$. Note
that $(H,\mathcal{O}_H)\subset (G_1,\mathcal{O}_{G_1})$ and
$(G/G_1,\mathcal{O}_{G/G_1})\simeq (\pt, \bigwedge W)$. Denote by
$\Phi$ the natural $(G,\mathcal{O}_{G})$-equivariant morphism
$(G/H,\mathcal{O}_{G/H})\to (G/G_1,\mathcal{O}_{G/G_1})$. Note that
$\Phi^*$ is injective, hence $\bigwedge W \simeq
\Phi^*(H^0(\mathcal{O}_{G/G_1}))\subset H^0(\mathcal{O}_{G/H})$.
Since $(M,\gr \mathcal{O}_M)$ is split,
$H^0(\gr\mathcal{O}_{M})\simeq \bigwedge W$. It is easy to see that
$\dim H^0(\mathcal{O}_{G/H})\leq \dim H^0(\gr\mathcal{O}_{G/H})$. It
follows that $\mathcal{O}_{G/H}\simeq \bigwedge W$.$\Box$

\bigskip

\begin{center}
{\it 2.4 An application of Theorem \ref{funk1}}
\end{center}

\medskip

 Let $(M,\mathcal O_M)$, $(B,\mathcal O_B)$ and $(F,\mathcal
O_F)$ be complex supermanifolds. The supermanifold $(M,\mathcal
O_M)$ is called a {\it bundle with fiber $(F,\mathcal O_F)$, base
space $(B,\mathcal O_B)$ and with projection $p: (M,\mathcal O_M)
\to (B,\mathcal O_B)$} if the following condition holds:
 there is an open covering $\{U_i\}$ of the manifold $B$ and
isomorphisms $\psi_i: (p_1^{-1}(U_i),\mathcal O_M)\rightarrow
(U_i,\mathcal O_B)\times (F,\mathcal O_F)$ such that the following
diagram is commutative:
$$
\begin{CD}
(p_1^{-1}(U_i),\mathcal O_M)@>{\psi_i}>>(U_i, \mathcal O_B)\times
(F,\mathcal O_F)\\
@V{ p}VV @VV {\pr_1}V\\
(U_i, \mathcal O_B)@=(U_i, \mathcal O_B)
\end{CD},
$$
where $\pr_1$ is the projection onto the first factor.

Let $(M,\mathcal{O}_M)$ be a supermanifold, $M$ a connected manifold
and $\mathfrak{g}$ a complex (finite dimensional) Lie superalgebra.
An action of $\mathfrak{g}$ on $(M,\mathcal{O}_M)$ is an arbitrary
Lie algebra homomorphism $\varphi: \mathfrak{g} \to
\mathfrak{v}(M,\mathcal{O}_M)$. Assume that $(M,\mathcal{O}_M)$ is a
bundle with base $(B,\mathcal{O}_B)$ and projection map $p$. A
natural question is under what conditions the action of
$\mathfrak{g}$ on $(M,\mathcal{O}_M)$ induces an action of
$\mathfrak{g}$ on $(B,\mathcal{O}_B)$.

\medskip

\t\label{bash} {\it  Let $p: (M,\mathcal O_M)\to (B,\mathcal O_B)$
be the projection of a superbundle with fiber $(F,\mathcal O_F)$. If
$H^0(\mathcal O_F) \simeq\Bbb C$, then any action of a Lie
superalgebra is projectable with respect to $p$.}

\medskip

This theorem was proved in \cite{bash} in the case when
$p:(M,\mathcal O_M)=(B,\mathcal O_B)\times (F,\mathcal O_F)\to
(B,\mathcal O_B)$ is the natural projection.  Obviously, it can be
generalized  to bundles.

\medskip

\begin{center}
{\bf{\large 3. Holomorphic functions on classical flag
supermanifolds }}
\end{center}



\begin{center}
{\it 3.1 Classical flag supermanifolds  }
\end{center}

\medskip

 Yu.I. Manin \cite{Man} introduced four series of compact complex homogeneous
supermanifolds corresponding to the following four series of
classical linear complex Lie superalgebras:
\begin{enumerate}
  \item $\mathfrak{gl}_{m|n}(\mathbb{C})$, the general linear Lie
  superalgebra of the vector superspace $\mathbb{C}^{m|n}$;
  \item $\mathfrak{osp}_{m|n}(\mathbb{C})$, the orthosymplectic Lie
  superalgebra that annihilates a non-degenerate even symmetric
  bilinear form in $\mathbb{C}^{m|n}$, $n$ even;
  \item $\pi\mathfrak{sp}_{n|n}(\mathbb{C})$, the linear Lie
  superalgebra that annihilates a non-degenerate odd skew-symmetric
  bilinear form in $\mathbb{C}^{n|n}$;
  \item $\mathfrak{q}_{n|n}(\mathbb{C})$, the linear Lie
  superalgebra that commutes with an odd involution $\pi$ in
  $\mathbb{C}^{n|n}$.
\end{enumerate}
These supermanifolds are called {\it supermanifolds of flags} in
Case $1$, {\it supermanifolds of isotropic flags} in Cases $2$ and
$3$, and {\it supermanifolds of $\pi$-symmetric flags} in Case $4$.
We will call all of them {\it classical flag supermanifolds}. For
further reading, see also \cite{Man_Top,Penkov_Skornyakov,Penkov}.

Denote by $\bold F^{m}_{k}$ the usual manifold of flags of type
$k=(k_1,\ldots,k_r)$ in $\Bbb C^m$, where $0 \leq k_r \leq\dots \leq
k_1 \leq m$. Let us describe an atlas on $\bold F^{m}_{k}$.

Let $\Bbb C^m\supset W_1\supset\dots\supset W_r$ be a flag of type
$k_1,\ldots,k_r$. Choose a basis $B_s$ in each $W_s$. Assume that
$B_0 = (e_1,\dots,e_m)$ is the standard basis of $\Bbb C^m$ and put
$k_0 = m$. Then for any $s = 1,\ldots,r$ the matrix $X_s\in
\operatorname{Mat}_{k_{s-1},k_s}(\Bbb C)$ is defined in the
following way: the columns of $X_s$ are the coordinates of the
vectors from $B_s$ with respect to the basis $B_{s-1}$. Since $\rk
X_s = k_s$, the matrix $X_s$ contains a non-degenerate minor of size
$k_s$.

For each $s = 1,\ldots,r$, let us fix a $k_s$-tuple
$I_s\subset\{1,\ldots,k_{s-1}\}$. Put $I = (I_1,\ldots,I_r)$. Denote
by $U_I$ the set of flags $f$ from $\bold F^{m}_{k}$ satisfying the
following conditions: there exist bases $B_s$ such that $X_s$
contains the identity matrix of size $k_s$ in the lines with numbers
from $I_s$. It is easy to see that any flag from $U_I$ is uniquely
determined by those elements of $X_s$ that are not contained in the
identity matrix. Furthermore, any flag is contained in a certain
 $U_I$. The elements of $X_s$ that are not contained in the identity
 matrix are the coordinates of a flag
 from $U_I$ in the chart determined by $I$. Rename $X_{I_s}:= X_s$.
 Hence the local coordinates in $U_I$ are determined by
 $r$-tuple $(X_{1},\ldots,X_{r})$. If $J = (J_1,\ldots,J_s)$,
 where $J_s\subset\{1,\ldots,k_{s-1}\},\;|J_s| = k_s$, then the
 transition functions between the charts $U_I$ and $U_J$
 are given by:
$$
X_{J_1} = X_{I_1}C_{I_1J_1}^{-1}, \ \ X_{J_s} =
C_{I_{s-1}J_{s-1}}X_{I_s}C_{I_sJ_s}^{-1},\ \ s\ge 2,
$$
where $C_{I_1J_1}$ is the submatrix of $X_{I_1}$ formed by the lines
with numbers from $J_1$ and $C_{I_sJ_s}$, $s\ge 2$, is the submatrix
of $C_{I_{s-1}J_{s-1}}X_{I_s}$ formed by lines with numbers from
$J_s$.

Let us give an explicit description of classical flag supermanifolds
in terms of atlases and local coordinates (see, \cite{V_flag,
V_flag_Pi, V_diss}). (Note that in \cite{Man} such a description was
given only for super-grassmannians.) Let us take $m,n\in\Bbb N$ and
let $k=(k_1,\ldots,k_r)$ and $l=(l_1,\ldots,l_r)$ be two $r$-tuples
such that  $0\le k_r\le\ldots\le k_1\le m,\; 0\le l_r\ldots\le
l_1\le n$ è $0 < k_r +l_r <\ldots < k_1 + l_1 < m+n$. Let us define
the {\it supermanifold $\bold F^{m|n}_{k|l}$ of flags of type }
$(k|l)$ in the superspace  $V = \Bbb C^{m|n}$. The reduction of
$\bold F^{m|n}_{k|l}$ will be the product $\bold F^{m}_{k}
\times\bold F^{n}_{l}$ of two manifolds of flags of type $k$ and $l$
in $\Bbb C^m = V_{\bar 0}$ and $\Bbb C^n = V_{\bar 1}$.

For each $s = 1,\ldots,r$, let us fix $k_s$- and $l_s$-tuples of
numbers
 $I_{s\bar 0}\subset\{1,\ldots,k_{s-1}\}$ and $I_{s\bar
1}\subset\{1,\ldots,l_{s-1}\}$, where $k_0=m$, $l_0=n$. We put
$I_s=(I_{s\bar 0},I_{s\bar 1}),\; I = (I_1,\ldots,I_r)$. Our aim is
now to construct a superdomain $\mathcal{W}_I$. To each $I_s$ assign
a matrix of size $(k_{s-1} + l_{s-1})\times (k_s + l_s)$
\begin{equation}
\label{Z_I_s}
\begin{split}
Z_{I_s}=\left(
\begin{array}{cc}
X_s & \Xi_s\\
\H_s & Y_s \end{array} \right), \ \ s=1,\dots,r.
\end{split}
\end{equation}
 Suppose that the identity matrix $E_{k_s+l_s}$ is contained in the
lines of $Z_{I_s}$ with numbers $i\in I_{s\bar 0}$ and $k_{s-1} +
j,\; j\in I_{s\bar 1}$. Here $X_s\in \operatorname{Mat}_{k_{s-1},
k_{s}}(\Bbb C),\; Y_s\in \operatorname{Mat}_{l_{s-1}, l_{s}}(\Bbb
C)$, where $\operatorname{Mat}_{a, b}(\Bbb C)$ is the space of
matrices of size $a\times b$ over $\Bbb C$. By definition, the
entries of $X_s$ and $Y_s$, $s=1,\dots,r$, that are not contained in
the identity matrix form the even coordinate system of
$\mathcal{W}_I$. The non-zero entries of $\Xi_s$ and $\H_s$ form the
odd coordinate system of $\mathcal{W}_I$.

Thus we have defined a set of superdomains on $\bold F^{m}_{k}
\times\bold F^{n}_{l}$ indexed by $I$. Note that the reductions of
these superdomains cover $\bold F^{m}_{k} \times\bold F^{n}_{l}$.
The local coordinates of each superdomain are determined by the
$r$-tuple of matrices $(Z_{I_1},\ldots,Z_{I_r})$. Let us define the
transition functions between two superdomains corresponding to $I =
(I_s)$ and $J = (J_s)$ by the following formulas:
\begin{equation}
\label{perehod}
\begin{split}
Z_{J_1} = Z_{I_1}C_{I_1J_1}^{-1}, \ \ Z_{J_s} =
C_{I_{s-1}J_{s-1}}Z_{I_s}C_{I_sJ_s}^{-1},\ \ s\ge 2,
\end{split}
\end{equation}
where $C_{I_1J_1}$ is the submatrix of $Z_{I_1}$ that consists of
the lines with numbers from $J_1$, and $C_{I_sJ_s}$, $s\ge 2$, is
the submatrix of $C_{I_{s-1}J_{s-1}}Z_{I_s}$ that consists of the
lines with numbers from $J_s$. Gluing the superdomains $\mathcal
W_I$, we define the {\it supermanifold of flags} $\bold
F^{m|n}_{k|l}$. In the case $r = 1$, this supermanifold is  called a
{\it super-grassmannian}. In the literature the notation $\bold
{Gr}_{m|n,k_1|l_1}$ is sometimes used.

The supermanifold $\bold F^{m|n}_{k|l}$ is $\GL_{m|n}(\mathbb
C)$-homogeneous. The action can be given by
\begin{equation}
\label{dey}
\begin{split}
(L,(Z_{I_1},\ldots,Z_{I_r}))\mapsto (\tilde Z_{J_1},\ldots,\hat
Z_{J_r}),\\
\tilde Z_{J_1} = LZ_{I_1}C_1^{-1},\,\, \tilde Z_{J_s} =
C_{s-1}Z_{I_s}C_s^{-1}.
\end{split}
\end{equation}
Here $L$ is a coordinate matrix of $\GL_{m|n}(\mathbb C)$, $C_1$ is
the invertible submatrix of $LZ_{I_1}$ that consists of the lines
with numbers from $J_1$, à $C_s,\; s\ge 2$, is the invertible
submatrix of $C_{s-1}Z_{I_s}$ that consists of the lines with
numbers from $J_s$.

Let $\mathfrak{g}$ be one of the classical Lie superalgebras
described in $3.1$. Denote by $\mathbf F_{k|l}(\mathfrak {g})$ the
flag supermanifold of type $(k|l)$ corresponding to $\mathfrak{g}$.
We will also write $\mathbf F_{k|l}(\mathfrak {gl}_{m|n}(\mathbb
C))= \bold F_{k|l}^{m|n}$. Let us describe
 $\bold F_{k|l}(\mathfrak {g})$ for $\mathfrak {g}=\mathfrak{osp}_{m|n}(\mathbb{C})$,
$\pi\mathfrak{sp}_{n|n}(\mathbb{C})$ or
$\mathfrak{q}_{n|n}(\mathbb{C})$ in coordinates.

The subsupermanifold $\mathbf F_{k|l}(\mathfrak {osp}_{m|2n}(\mathbb
C))$ of $\bold F^{m|2n}_{k|l}$ is given in coordinates (\ref{Z_I_s})
by the following equations:
\begin{equation}
\begin{split}
\left(
\begin{array}{cc}
X_1 & \Xi_1\\
\H_1 & Y_1 \end{array} \right)^{ST} \Gamma \left(
\begin{array}{cc}
X_1 & \Xi_1\\
\H_1 & Y_1 \end{array} \right)=0,
\end{split}
\end{equation}
where
\begin{equation}
\label{forma even} \Gamma=\left(
  \begin{array}{cccc}
    0 & E_s & 0 & 0 \\
    E_s & 0 & 0 & 0 \\
    0 & 0 & 0 & E_n \\
    0 & 0 & -E_n & 0 \\
  \end{array}
\right), \quad \Gamma=\left(
  \begin{array}{ccccc}
  1&0 & 0 & 0 & 0 \\
   0& 0 & E_s & 0 & 0 \\
    0&E_s & 0 & 0 & 0 \\
    0&0 & 0 & 0 & E_n \\
    0&0 & 0 & -E_n & 0 \\
  \end{array}
\right),
\end{equation}
$m=2s$ or $m=2s+1$ and
$$
\left(
\begin{array}{cc}
X & \Xi\\
\H & Y \end{array} \right)^{ST}= \left(
\begin{array}{cc}
X^T & \H^T\\
-\Xi^T & Y^T \end{array} \right)
$$
 is the super-transposition.

The subsupermanifold $\mathbf F_{k|l}(\pi\mathfrak{sp}_{n}(\mathbb
C))$ of $\bold F^{n|n}_{k|l}$ is given in coordinates (\ref{Z_I_s})
by the following equations:
\begin{equation}\label{pisp_uravnen}
\begin{split}
\left(
\begin{array}{cc}
X_1 & \Xi_1\\
\H_1 & Y_1 \end{array} \right)^{ST} \Upsilon \left(
\begin{array}{cc}
X_1 & \Xi_1\\
\H_1 & Y_1 \end{array} \right)=0,
\end{split}
\end{equation}
where
\begin{equation}\label{forma_pisp}
\Upsilon=\left(
  \begin{array}{cc}
    0 & E_n \\
   - E_n & 0 \\
  \end{array}
\right)
\end{equation}

The subsupermanifold $\mathbf F_{k|l}(\mathfrak {q}_{n}(\mathbb C))$
of $\bold F^{n|n}_{k|k}$ is given in coordinates (\ref{Z_I_s}) by
 $X_s=Y_s$, $\Xi_s=\H_s$, $s=1,\ldots, r$.

In \cite{Man} the action of $(G,\mathcal O_G)=\OSp_{m|2n}(\Bbb C)$,
$\Pi\!\Sp_n(\Bbb C)$ or $\Q_n(\Bbb C)$ on $\bold F_{k|l}(\mathfrak
{g})$ was defined. In our coordinates this action is given by
$(\ref{dey})$, where we assume that $L$ is a cordinate matrix of
$(G,\mathcal O_G)$.

\medskip

\begin{center}
{\it 3.2 Holomorphic functions on classical flag supermanifolds  }
\end{center}

\medskip

To compute the algebra of holomorphic functions on
$\mathbf{F}_{k|l}(\mathfrak{g})$ using Theorem \ref{funk1} we need
to know the Lie superalgebra $\mathfrak{p}_{\mathfrak{g}}$ of the
stabilizer of a point $x\in (\mathbf{F}_{k|l}(\mathfrak{g}))_{\red}$
for the action (\ref{dey}). Such stabilizers  are also called
parabolic subsupergroups of $(G,\mathcal O_G)=\OSp_{m|2n}(\Bbb C)$,
$\Pi\!\Sp_n(\Bbb C)$ or $\Q_n(\Bbb C)$ (see \cite{Man_Top, Penkov}).
We follow the approach of A.L.~Onishchik  \cite{Ivan_Onishch}. We
will use the following lemma.

\medskip

\lem\label{gamma is inecive} {\it In the conditions of Theorem
\ref{funk1}, assume that $\mathfrak{g}_{\bar 1}$ is a completely
reducible $G$-module. Let $W\subset \Im \gamma^*$ be a non-trivial
$G$-module. Denote by $X$ the following $G$-module
$$
X:=\{ v\in \mathfrak{g}_{\bar 1}\,\,\mid\,\, W(v)= \{0\}\}
$$
and by $Y$ a complement to $X$ in $\mathfrak{g}_{\bar 1}$. Then
$\gamma|_Y$ is injective.

 }

\medskip

\noindent{\it Proof.} Assume that $\gamma(v)=0$ for some $v\in
Y\backslash \{0\}$. Then there is an $f\in W$ such that $f(v)\ne 0$.
(Otherwise $W(v)=\{0\}$ and $v\in X$.) Since there exists an $l\in
(\mathfrak{g}_{\bar 1}/\mathfrak{h}_{\bar 1})^*$ such that
$\gamma^*(l)=f$, we arrive at a contradiction.$\Box$

\bigskip

\noindent{\bf Case $\mathfrak{g}=\mathfrak{gl}_{m|n}(\mathbb{C})$.}
Let $e_1,\ldots,e_m,f_1,\ldots, f_n$ be the standard basis of
$\mathbb{C}^{m|n}$. Consider the superdomain $Z_{I}$ in $\bold
F_{k|l}(\mathfrak{gl}_{m|n}(\mathbb{C}))$ corresponding to $I_{s\bar
0} = (1,\ldots,k_s)$, $I_{s\bar 1} = (1,\ldots,l_s)$. Denote by $x$
the origin of $Z_{I}$. It is easy to see that $x=(V_1,\ldots,V_r)$,
where $V_i=\langle e_1,\ldots,e_{k_i}\rangle \oplus \langle
f_1,\ldots,f_{l_i}\rangle$. Denote by
$\mathfrak{p}(x)_{\mathfrak{gl}}$ the Lie superalgebra of the
stabilizer of $x$ for the action (\ref{dey}) of
$\GL_{m|n}(\mathbb{C})$. It is easy to see that
$$
\mathfrak{p}(x)_{\mathfrak{gl}} =\{X\in
\mathfrak{gl}_{m|n}(\mathbb{C})\,\mid \, X(V_i)\subset V_i\}.
$$

The Lie superalgebra $\mathfrak{p}(x)_{\mathfrak{gl}}$ admits
another description in terms of root systems, see
\cite{Ivan_Onishch}, which we are going to describe now. Let us take
a Cartan subalgebra $\mathfrak{t}\subset
\mathfrak{gl}_{m|n}(\mathbb{C})_{\bar 0}$ in the following form
$$
\begin{array}{l}
\operatorname{diag}(x_1,\ldots,x_m,y_1,\ldots,y_n).
\end{array}
$$
The corresponding root system $\Delta=\Delta_{\bar 0} \cup
\Delta_{\bar 1}$ is given by
$$
\begin{array}{l}
\Delta_{\bar 0}=\{x_i- x_j, y_i- y_j\, \, \mid\, i\ne j \},\,\,\,
\Delta_{\bar 1}= \{x_i- y_j, \, y_i-x_j\}.
\end{array}
$$
Let us take an $m$-tuple $a=(a_1,\ldots,a_m)$ and an $n$-tuple
$b=(b_1,\ldots,b_n)$ of real numbers such that
$$
\begin{array}{c}
a_1=\cdots=a_{k_r}=b_1=\cdots=b_{l_r}> \cdots>
a_{k_{2}+1}=\cdots=a_{k_1}=\\
b_{l_2+1}=\cdots=b_{l_1}>
a_{k_{1}+1}=\cdots=a_{m}=b_{l_1+1}=\cdots=b_{n}.
\end{array}
$$
Then $(a,b)\in \mathfrak{t}(\mathbb{R})$. Let
\begin{equation}
\label{parabolic_gl} \mathfrak{p}(a,b)_{\mathfrak{gl}}=
\mathfrak{t}\oplus \bigoplus_{\alpha\in \Delta,\,\alpha(a,b)\geq 0}
\mathfrak{gl}_{m|n}(\mathbb{C})_{\alpha}.
\end{equation}
Note that $\mathfrak{p}(a,b)_{\mathfrak{gl}}$ depends only on the
numbers $k_i,\, l_i$, $i=1,\ldots,r$.
 From \cite{Ivan_Onishch},
Chapter $4$, $\S\,1$, Proposition $1$, it can be deduced that
$\mathfrak{p}(a,b)_{\mathfrak{gl}} =
\mathfrak{p}(x)_{\mathfrak{gl}}$. (This follows also from the direct
calculation.)



\medskip

\t\label{theor_funk gl} {\it If $ (k|l)\ne (m,\ldots, m,
k_{s+2},\ldots, k_r)| (l_1,\ldots, l_s,0,\ldots,0) $ and
$$
(k|l)\ne (k_1,\ldots, k_s,0,\ldots,0) |(n,\ldots, n, l_{s+2},\ldots,
l_r),
$$
then $H^0(\mathbf{F}_{k|l}(\mathfrak{gl}_{m|n}(\mathbb{C})))\simeq
\mathbb{C}$. Otherwise
$$
\mathbf{F}_{k|l}(\mathfrak{gl}_{m|n}(\mathbb{C})) \simeq
(\pt,\bigwedge(mn))\times (\mathbf{F}_{k}\times \mathbf{F}_{l})
$$
and  $H^0(\mathbf{F}_{k|l}(\mathfrak{gl}_{m|n}(\mathbb{C})))\simeq
\bigwedge(mn)$.

}

\medskip

\noindent{\it Proof.} The odd part
$\mathfrak{gl}_{m|n}(\mathbb{C})_{\bar 1}$ of the Lie superalgebra
$\mathfrak{gl}_{m|n}(\mathbb{C})$ for $m,n\geq 1$ is the direct sum
of two irreducible  $\mathfrak{gl}_{m|n}(\mathbb{C})_{\bar
0}$-submodules
$$
V_1=\left\{ \left(
  \begin{array}{cc}
    0 & A \\
    0& 0 \\
  \end{array}
\right), \,\,\, A\in \operatorname{Mat}_{m\times
n}(\mathbb{C})\right\},\,V_2= \left\{ \left(
  \begin{array}{cc}
    0 & 0 \\
    B& 0 \\
  \end{array}
\right), \, B\in \operatorname{Mat}_{n\times m}(\mathbb{C})
\right\},
$$
and this decomposition is unique.

If $ (k|l) = (m,\ldots, m, k_{s+2},\ldots, k_r)| (l_1,\ldots,
l_s,0,\ldots,0) $ then the Lie superalgebra
$\mathfrak{p}(x)_{\mathfrak{gl}}=\mathfrak{p}(a,b)_{\mathfrak{gl}}$
is determined by an $m$-tuple $a=(a_1,\ldots, a_m)$ and an $n$-tuple
$b=(b_1,\ldots, b_n)$ such that
$$
a_1\geq \cdots\geq a_m > b_1\geq \cdots\geq b_n.
$$
Hence $\mathfrak{p}(x)_{\mathfrak{gl}}\supset V_1$. Consider the
subalgebra $\mathfrak{g}' = \mathfrak{gl}_{m|n}(\mathbb{C})_{\bar
0}\oplus V_2$ in $\mathfrak{gl}_{m|n}(\mathbb{C})$. Denote by
$(G',\mathcal{O}_{G'})$ the subsupergroup of $\GL_{m|n}(\mathbb{C})$
 determined by the Harish-Chandra pair $(G,\mathfrak{g}')$.
It is easy to see that $(G',\mathcal{O}_{G'})$ acts on
$\mathbf{F}_{k|l}(\mathfrak{gl}_{m|n}(\mathbb{C}))$ transitively for
this $(k|l)$ and the Lie superalgebra
$\mathfrak{p}'_{\mathfrak{gl}}$ of the stabilizer of $x$ for this
action is $\mathfrak{p}(x)_{\mathfrak{gl}}\cap
\mathfrak{g}'=(\mathfrak{p}(x)_{\mathfrak{gl}})_{\bar 0}$. It
follows that the stabilizer is a usual Lie group. We get that
$\mathbf{F}_{k|l}(\mathfrak{gl}_{m|n}(\mathbb{C}))$ is split and the
structure sheaf of
$\mathbf{F}_{k|l}(\mathfrak{gl}_{m|n}(\mathbb{C}))$ is isomorphic to
$\mathcal{F}_{M}\otimes \bigwedge V_2^*$, where $M =
(\mathbf{F}_{k|l}(\mathfrak{gl}_{m|n}(\mathbb{C})))_{\red}$, see
Example \ref{example2}. In particular,
$$
H^0(\mathbf{F}_{k|l}(\mathfrak{gl}_{m|n}(\mathbb{C}))) \simeq
\bigwedge V_2^* \simeq \bigwedge(mn).
$$

In the case $ (k|l) = (k_1,\ldots, k_s,0,\ldots,0) |(n,\ldots, n,
l_{s+2},\ldots, l_r) $ the proof is similar.

Assume that
$$
\begin{array}{l}
 (k|l) \ne (m,\ldots, m, k_{s+2},\ldots, k_r)| (l_1,\ldots,
l_s,0,\ldots,0)\,\, \text{and} \\
 (k|l) \ne  (k_1,\ldots,
k_s,0,\ldots,0) |(n,\ldots, n, l_{s+2},\ldots, l_r).
\end{array}
$$
 Then there
are $i,j$ such that $a_i=b_j$ or there are $i_1,j_1$ and $i_2,j_2$
such that $a_{i_1}>b_{j_1}$ and $a_{i_2} < b_{j_2}$. In the first
case, $\Ker \gamma$ contains the subspace $\langle E_{i,m+j},
E_{m+j,i}\rangle$. In the second case, $\Ker \gamma$ contains the
subspace $\langle E_{i_1,m+j_1}, E_{m+j_2,i_2}\rangle$. Thus,
$\gamma|_{V_k}$, $k=1,2$, cannot be injective.
 By Lemma
\ref{gamma is inecive} and Theorem \ref{funk1}, we get that
$H^0(\mathbf{F}_{k|l}(\mathfrak{gl}_{m|n}(\mathbb{C})))\simeq
\mathbb{C}$.$\Box$

\bigskip

\noindent{\bf Case $\mathfrak{g}=\mathfrak{osp}_{m|n}(\mathbb{C})$.}
Since the manifold
$\mathbf{F}_{k|l}(\mathfrak{osp}_{m|n}(\mathbb{C}))_{\red}$ consists
of isotropic flags, it follows that $k_1\leq p:=[\frac{m}{2}]$ and
$l_1\leq q:=\frac{n}{2}$. Let us fix a basis of $\mathbb{C}^{m|n}$
such that the matrix of the corresponding non-degenerate even
symmetric bilinear form has the matrix $\Gamma$ given by (\ref{forma
even}) and denote its elements as follows:
$$
\begin{array}{l}
e_1,\ldots,e_{2p}, f_1,\ldots f_n, \,\,\text{if $m$ is even},\\
e_0,\ldots,e_{2p}, f_1,\ldots f_n, \,\,\text{if $m$ is odd}.
\end{array}
$$
Consider the superdomain $Z_{I}$ in $\bold F^{m|n}_{k|l}$
corresponding to $I_{s\bar 0} = (1,\ldots,k_s)$ and $I_{s\bar 1} =
(1,\ldots,l_s)$. Denote by $x$ the origin of $Z_{I}$. It is easy to
see that $x=(V_1,\ldots,V_r)$, where $V_i=\langle
e_1,\ldots,e_{k_i}\rangle \oplus \langle f_1,\ldots,f_{l_i}\rangle$,
and $x$ is isotropic. Denote by $\mathfrak{p}(x)_{\mathfrak{osp}}$
the Lie superalgebra of the stabilizer of $x$ for the action
(\ref{dey}) of $\OSp_{m|n}(\mathbb{C})$. It is easy to see that
$$
\mathfrak{p}(x)_{\mathfrak{osp}} =\{X\in
\mathfrak{osp}_{m|n}(\mathbb{C})\,\mid \, X(V_i)\subset V_i\}.
$$

The Lie superalgebra $\mathfrak{osp}_{m|n}(\mathbb{C})$ has the
following forms for $m=2p+1$ or $m=2p$, respectively:
$$
\begin{array}{c}
\left(
  \begin{array}{ccccc}
    0 & -v^t & -u^t & w & w_1 \\
    u & A & B & U & U_1 \\
    v & C & -A^t & W & W_1\\
    w_1^t & W_1^t & U_1^t & Y & Z \\
    -w^t & -W^t & -U^t & T & -Y^t \\
  \end{array}
\right), \,\,\,\, \left(
  \begin{array}{cccc}
    A & B & U & U_1 \\
     C & -A^t & W & W_1\\
     W_1^t & U_1^t & Y & Z \\
     -W^t & -U^t & T & -Y^t \\
  \end{array}
\right),\\
 B^t=-B,\, C^t=-C,\,\, Z^t=Z,\,\, T^t=T.
\end{array}
$$
Here $Y,Z,T$ are square matrices of order $q$, $A,B,C$ are square
matrices of order $p$, $U,U_1,V,V_1$ are $p\times q$-matrices, $u,v$
are columns of height $p$, and $w,w_1$ are rows of length $q$. As a
Cartan subalgebra $\mathfrak{t}\subset
\mathfrak{osp}_{m|n}(\mathbb{C})_{\bar 0}$, one takes that of all
diagonal matrices
$$
\begin{array}{c}
\operatorname{diag}(x_1,\ldots,x_p,-x_1,\ldots,-x_p,
y_1,\ldots,y_q,-y_1,\ldots,-y_q), \,\,\text{for}\,\, m=2p,\\
\operatorname{diag}(0,x_1,\ldots,x_p,-x_1,\ldots,-x_p,
y_1,\ldots,y_q,-y_1,\ldots,-y_q), \,\,\text{for}\,\, m=2p+1.
\end{array}
$$
The corresponding root system $\Delta=\Delta_{\bar 0} \cup
\Delta_{\bar 1}$ is given by
$$
\begin{array}{l}
\Delta_{\bar 0}=\left\{
                  \begin{array}{ll}
                    \{\pm x_i\pm x_j,\,\pm y_i\pm y_j, \, \pm 2y_i\, \mid\, i\ne j \} & \hbox{for $m=2p$,} \\
                    \{\pm x_i\pm x_j,\, \pm x_i,\,\pm y_i\pm y_j, \, \pm 2y_i\, \mid\, i\ne j \} & \hbox{for $m=2p+1$,}
                  \end{array}
                \right.\\
\Delta_{\bar 1}=\left\{
                  \begin{array}{ll}
                    \{\pm x_i\pm y_j \} & \hbox{for $m=2p$,} \\
                    \{\pm x_i\pm y_j,\, \pm y_i \} & \hbox{for $m=2p+1$.}
                  \end{array}
                \right.
\end{array}
$$

 Let us take a $p$-tuple $a=(a_1,\ldots,a_p)$ and a
$q$-tuple $b=(b_1,\ldots,b_q)$ of real numbers such that
$$
\begin{array}{c}
a_1=\cdots=a_{k_r}=b_1=\cdots=b_{l_r}> \cdots>
a_{k_{2}+1}=\cdots=a_{k_1}=\\
b_{l_2+1}=\cdots=b_{l_1}>
a_{k_{1}+1}=\cdots=a_{p}=b_{l_1+1}=\cdots=b_{p}=0.
\end{array}
$$
Then $(a,-a,b,-b)\in \mathfrak{t}(\mathbb{R})$. Let
\begin{equation}
\label{parabolic_osp} \mathfrak{p}(a,b)_{\mathfrak{osp}}=
\mathfrak{t}\oplus \bigoplus_{\alpha\in
\Delta,\,\alpha(a,-a,b,-b)\geq 0}
\mathfrak{osp}_{m|n}(\mathbb{C})_{\alpha} .
\end{equation}
Note that $\mathfrak{p}(a,b)_{\mathfrak{osp}}$ depends only on the
numbers $k_i,\, l_i$, $i=1,\ldots,r$. In \cite{Ivan_Onishch},
Chapter $4$, $\S\,2$, Proposition $2$, it was shown that
$\mathfrak{p}(a,b)_{\mathfrak{osp}}=\mathfrak{p}(x)_{\mathfrak{osp}}$
if $m\geq 1$, $n\geq 2$.


\medskip

\t\label{theor_funk osp} {\it Assume that $m\geq 1$, $n\geq 2$.
 If $m$ is odd or $m$ is even and $m>2$, then
$H^0(\mathbf{F}_{k|l}(\mathfrak{osp}_{2|n}(\mathbb{C}))) \simeq
\mathbb{C}$.

Suppose that $m=2$. If $k|l\ne (1,\ldots,1|l_1,\ldots, l_{r-1},0)$,
then $$ H^0(\mathbf{F}_{k|l}(\mathfrak{osp}_{2|n}(\mathbb{C})))
\simeq \mathbb{C}.
$$

Suppose that $m=2$ and $k|l= (1,\ldots,1|l_1,\ldots, l_{r-1},0)$,
then
$$
\mathbf{F}_{k|l}(\mathfrak{osp}_{2|n}(\mathbb{C}))\simeq
(\pt,\bigwedge(2q))\times M,
$$
where
$M=(\mathbf{F}_{k|l}(\mathfrak{osp}_{2|n}(\mathbb{C})))_{\red}$. In
particular,
$H^0(\mathbf{F}_{k|l}(\mathfrak{osp}_{2|n}(\mathbb{C})))\simeq
\bigwedge(2q)$.

}

\medskip

\noindent{\it Proof.} Assume that $m$ is odd or $m$ is even and
$m>2$, then $\mathfrak{osp}_{m|n}(\mathbb{C})_{\bar 1}$ is an
irreducible $\mathfrak{osp}_{m|n}(\mathbb{C})_{\bar 0}$-module.
Hence by Lemma \ref{gamma is inecive} it is sufficient to check that
$(\mathfrak{p}(a,b)_{\mathfrak{osp}})_{\bar 1}\ne \{0\}$. Since
$a_i,b_j\geq 0$, we get that
$\mathfrak{osp}_{m|n}(\mathbb{C})_{x_i+y_j}\subset
(\mathfrak{p}(a,b)_{\mathfrak{osp}})_{\bar 1}$ and
$\mathfrak{osp}_{m|n}(\mathbb{C})_{y_j}\subset
(\mathfrak{p}(a,b)_{\mathfrak{osp}})_{\bar 1}$.

Assume that $m=2$. In this case
$\mathfrak{osp}_{m|n}(\mathbb{C})_{\bar 1}$ is a direct sum of two
irreducible $\mathfrak{osp}_{m|n}(\mathbb{C})_{\bar 0}$-modules:
$$
V_1= \left(
       \begin{array}{cccc}
         0 & 0 & 0 & 0 \\
         0 & 0 & W & W_1 \\
         W_1^t & 0 & 0 & 0 \\
         -W^t & 0 & 0 & 0 \\
       \end{array}
     \right),\quad
V_2= \left(
       \begin{array}{cccc}
         0 & 0 & U & U_1 \\
         0 & 0 & 0 & 0 \\
         0 & U_1^y & 0 & 0 \\
         0& -U^t & 0 & 0 \\
       \end{array}
     \right),
$$
and this decomposition is unique. Assume that $k|l=
(1,\ldots,1|l_1,\ldots, l_{r-1},0)$. Then the Lie superalgebra
$\mathfrak{p}(x)_{\mathfrak{osp}}=\mathfrak{p}(a,b)_{\mathfrak{osp}}$
is determined by a $1$-tuple $a=(a_1)$ and a $q$-tuple
$b=(b_1,\ldots, b_q)$ such that
$$
a_1 > b_1\geq \cdots\geq b_q\geq 0.
$$
In this case  $(\mathfrak{p}(a,b)_{\mathfrak{osp}})_{\bar 1}=V_2$.
Consider the subalgebra $\mathfrak{g}' =
\mathfrak{osp}_{m|n}(\mathbb{C})_{\bar 0}\oplus V_1$ in
$\mathfrak{osp}_{m|n}(\mathbb{C})$. Denote by
$(G',\mathcal{O}_{G'})$ the subsupergroup of
$\OSp_{m|n}(\mathbb{C})$ determined by the Harish-Chandra pair
$(G,\mathfrak{g}')$. It is clear that $(G',\mathcal{O}_{G'})$ acts
on $\mathbf{F}_{k|l}(\mathfrak{osp}_{m|n}(\mathbb{C}))$ transitively
for this $(k|l)$ and the Lie superalgebra
$\mathfrak{p}'_{\mathfrak{osp}}$ of the stabilizer of $x$ for this
action is $\mathfrak{p}(x)_{\mathfrak{osp}}\cap
\mathfrak{g}'=(\mathfrak{p}(x)_{osp})_{\bar 0}$. It follows that the
stabilizer is a usual Lie group. We get that
$\mathbf{F}_{k|l}(\mathfrak{osp}_{m|n}(\mathbb{C}))$ is split and
the structure sheaf of
$\mathbf{F}_{k|l}(\mathfrak{osp}_{m|n}(\mathbb{C}))$ is isomorphic
to $\mathcal{F}_{M}\otimes \bigwedge V_1^*$, where $M=
(\mathbf{F}_{k|l}(\mathfrak{osp}_{m|n}(\mathbb{C})))_{\red}$, see
Example \ref{example2}. In particular,
$$
H^0(\mathbf{F}_{k|l}(\mathfrak{osp}_{m|n}(\mathbb{C}))) \simeq
\bigwedge V_1^* \simeq \bigwedge(mn).
$$

Assume that $k|l\ne  (1,\ldots,1|l_1,\ldots, l_{r-1},0)$. Then we
have the following possibilities:
\begin{itemize}
  \item there exists $j$ such that $a_1=b_j$,
  \item there exist $i,\,j$ such
that $a_{1}>b_{i}$ and $a_{1} < b_{j}$,
  \item $b_1\geq \cdots\geq b_n >
a_1$
\end{itemize}
 In the first case, $\gamma|_{V_k}$, $k=1,2$, cannot be
injective because $\mathfrak{osp}_{m|n}(\mathbb{C})_{x_1-y_i}\subset
V_2$, $\mathfrak{osp}_{m|n}(\mathbb{C})_{-x_1+y_i}\subset V_1$ and
$\mathfrak{osp}_{m|n}(\mathbb{C})_{x_1-y_i}\oplus
\mathfrak{osp}_{m|n}(\mathbb{C})_{-x_1+y_i}\subset
(\mathfrak{p}(x)_{\mathfrak{osp}})_{\bar 1}$. In the second case,
$\gamma|_{V_k}$, $k=1,2$, cannot be injective because again
$$
\mathfrak{osp}_{m|n}(\mathbb{C})_{x_1-y_i}\oplus
\mathfrak{osp}_{m|n}(\mathbb{C})_{-x_1+y_j}\subset
(\mathfrak{p}(x)_{\mathfrak{osp}})_{\bar 1}.$$
 In the third case,
$\gamma|_{V_k}$, $k=1,2$, cannot be injective because
$$
\mathfrak{osp}_{m|n}(\mathbb{C})_{x_1+y_i}\oplus
\mathfrak{osp}_{m|n}(\mathbb{C})_{-x_1+y_j}\subset
(\mathfrak{p}(x)_{\mathfrak{osp}})_{\bar 1}. \Box
$$

\bigskip

\noindent{\bf Case
$\mathfrak{g}=\pi\mathfrak{sp}_{n|n}(\mathbb{C})$.} The manifold
$\mathbf{F}_{k|l}(\pi\mathfrak{sp}_{n|n}(\mathbb{C}))_{\red}$
consists of isotropic flags, so $k_1+l_1\leq n$. Let us fix a basis
of $\mathbb{C}^{m|n}$ such that the matrix of the corresponding
non-degenerate odd symmetric bilinear form has the matrix
$\Upsilon$, see (\ref{forma_pisp}), and denote its elements as
follows:
$$
\begin{array}{l}
e_1,\ldots,e_{n}, f_1,\ldots f_n.
\end{array}
$$
Consider the superdomain $Z_{I}$ in $\bold F^{n|n}_{k|l}$
corresponding to $I_{s\bar 0} = (1,\ldots,k_s)$ and $I_{s\bar 1} =
(n-l_s+1,\ldots,n)$. Denote by $x$ the origin of $Z_{I}$. It is easy
to see that $x=(V_1,\ldots,V_r)$, where $V_i=\langle
e_1,\ldots,e_{k_i}\rangle \oplus \langle
f_{n-l_i+1},\ldots,f_{n}\rangle$, and $x$ is isotropic. Denote by
$\mathfrak{p}(x)_{\pi\mathfrak{sp}}$ the Lie superalgebra of the
stabilizer of $x$ for the action (\ref{dey}) of
$\Pi\!\Sp_{n|n}(\mathbb{C})$ on $\bold F_{k|l}
(\pi\mathfrak{sp}_{n|n}(\mathbb{C}))$. Again
$$
\mathfrak{p}(x)_{\pi\mathfrak{sp}} =\{X\in
\pi\mathfrak{sp}_{n|n}(\mathbb{C})\,\mid \, X(V_i)\subset V_i\}.
$$

The Lie superalgebra $\pi\mathfrak{sp}_{n|n}(\mathbb{C})$ has the
following form:
$$
\left(
  \begin{array}{cc}
    X & Y \\
    Z & -X^t \\
  \end{array}
\right),\,\,\, X,Y,Z \in \mathfrak{gl}_n(\mathbb{C}),\,\, Y=-Y^t,\,
Z=Z^t.
$$
As a Cartan subalgebra $\mathfrak{t}\subset
\pi\mathfrak{sp}_{n|n}(\mathbb{C})_{\bar 0}$, one takes that of all
diagonal matrices
$$
\begin{array}{l}
\operatorname{diag}(x_1,\ldots,x_n,-x_1,\ldots,-x_n).
\end{array}
$$
The corresponding root system $\Delta=\Delta_{\bar 0} \cup
\Delta_{\bar 1}$ is given by
$$
\begin{array}{l}
\Delta_{\bar 0}=\{x_i- x_j\, \mid\, i\ne j \},\,\,\, \Delta_{\bar
1}= \{x_i+ x_j\, \mid\, i< j \} \cup \{-x_i- x_j\, \mid\, i\leq j
\}.
\end{array}
$$

 Let us take an $n$-tuple $a=(a_1,\ldots,a_n)$ of real numbers such that
$$
\begin{array}{c}
|a_1|=\cdots=|a_{k_r}|=|a_{n-l_r+1}|=\cdots=|a_{n}|
>\\
|a_{k_r+1}|=\cdots=|a_{k_{r-1}}|=|a_{n-l_{r-1}+1}|=\cdots=|a_{n-l_r}|>\cdots>
\\
|a_{k_2}+1|=\cdots=|a_{k_1}|=|a_{n-l_1+1}|=\cdots=|a_{n-l_2}|>\\
a_{k_1+1} = \cdots =a_{n-l_1}=0,\\
a_i>0,\,\,\text{if} \,\, i\in \{1,\ldots,k_1\},\\
a_i<0,\,\,\text{if} \,\, i\in \{n-l_1+1,\ldots,n\}.
\end{array}
$$
Then $(a,-a)\in \mathfrak{t}(\mathbb{R})$. Let
\begin{equation}
\label{parabolic_psp} \mathfrak{p}(a)_{\pi\mathfrak{sp}}=
\mathfrak{t}\oplus \bigoplus_{\alpha\in \Delta,\,\alpha(a,-a)\geq 0}
\pi\mathfrak{sp}_{n|n}(\mathbb{C})_{\alpha} .
\end{equation}
Note that $\mathfrak{p}(a)_{\pi\mathfrak{sp}}$ depends only on the
numbers $k_i,\, l_i$, $i=1,\ldots,r$. In \cite{Ivan_Onishch},
Chapter $4$, $\S\,3$, Proposition $3$, it was shown that
$\mathfrak{p}(a)_{\pi\mathfrak{sp}}=\mathfrak{p}(x)_{\pi\mathfrak{sp}}$
if $n\geq 2$.


\medskip

\t\label{theor_funk pisp} {\it $1.$ Assume that $n\geq 2$.

If $k|l=(n,k_2,\ldots, k_r|0,\ldots, 0)$, then $H^0(\bold F_{k|l}
(\pi\mathfrak{sp}_{n|n}(\mathbb{C}))) \simeq \bigwedge ((n+1)n/2)$.

If $k|l=(0,\ldots, 0|n,l_2,\ldots, l_r)$ or $(0,\ldots,
0|n-1,l_2,\ldots, l_r)$ or $(1,0,\ldots, 0|n-1,n-1,l_3,\ldots,
l_r)$, then $H^0(\bold F_{k|l} (\pi\mathfrak{sp}_{n|n}(\mathbb{C})))
\simeq \bigwedge ((n-1)n/2)$.

For other $k|l$ we have $H^0(\bold F_{k|l}
(\pi\mathfrak{sp}_{n|n}(\mathbb{C}))) \simeq \mathbb{C}$.

$2.$ Assume that $n=1$, then $k|l=(1|0)$ or $(0|1)$. We have
$$
\begin{array}{c}
\bold F_{1|0} (\pi\mathfrak{sp}_{1|1}(\mathbb{C}))\simeq (\pt,
\bigwedge (1)) \,\,\,\text{and}\,\,\, H^0(\bold F_{1|0}
(\pi\mathfrak{sp}_{1|1}(\mathbb{C}))) \simeq
\mathbb{C} \oplus \mathbb{C},\\
\bold F_{0|1} (\pi\mathfrak{sp}_{1|1}(\mathbb{C}))\simeq
(\pt,\mathbb{C}) \,\,\,\text{and}\,\,\,  H^0(\bold F_{0|1}
(\pi\mathfrak{sp}_{1|1}(\mathbb{C}))) \simeq \mathbb{C}.
\end{array}
$$

}

\medskip

\noindent{\it Proof.} Suppose that $n\geq 2$. Then
$\pi\mathfrak{sp}_{n|n}(\mathbb{C})_{\bar 1}$ is a direct sum of two
irreducible $\pi\mathfrak{sp}_{n|n}(\mathbb{C})_{\bar 0}$-modules
$$
V_1= \left\{\left(
  \begin{array}{cc}
    0 & 0 \\
    Z & 0 \\
  \end{array}
\right) \right\}, \quad V_2= \left\{ \left(
  \begin{array}{cc}
    0 & Y \\
    0 & 0 \\
  \end{array}
\right) \right\},
$$
and this decomposition is unique.

Assume that $\gamma|_{V_1}$ is injective. Then $a_i+a_j> 0$ for all
$i\leq j$ and $k|l= (n,k_2,\ldots, k_r|0,\ldots, 0)$. Hence,
$(\mathfrak{p}(a)_{\pi\mathfrak{sp}})_{\bar 1}=V_2$. Consider the
subsuperalgebra $\mathfrak{g}' =
\pi\mathfrak{sp}_{n|n}(\mathbb{C})_{\bar 0}\oplus V_1$ in
$\pi\mathfrak{sp}_{n|n}(\mathbb{C})$. Denote by
$(G',\mathcal{O}_{G'})$ the subsupergroup of
$\Pi\Sp_{m|n}(\mathbb{C})$ determined by the Harish-Chandra pair
$(G,\mathfrak{g}')$. It is clear that $(G',\mathcal{O}_{G'})$ acts
on $\bold F_{k|l} (\pi\mathfrak{sp}_{n|n}(\mathbb{C}))$ transitively
and the Lie superalgebra of the stabilizer of $x$ for this action is
$\mathfrak{p}(x)_{\pi\mathfrak{sp}_{n|n}(\mathbb{C})}\cap
\mathfrak{g}'=
(\mathfrak{p}(x)_{\pi\mathfrak{sp}_{n|n}(\mathbb{C})})_{\bar 0}$.
Since $(\mathfrak{p}(x)_{\pi\mathfrak{sp}_{n|n}(\mathbb{C})})_{\bar
0}$ is a Lie algebra, we get that $\bold F_{k|l}
(\pi\mathfrak{sp}_{n|n}(\mathbb{C}))$ is split, see Example
\ref{example2}, and the structure sheaf of $\bold F_{k|l}
(\pi\mathfrak{sp}_{n|n}(\mathbb{C}))$ is isomorphic to
$\mathcal{F}_M\otimes \bigwedge V_1^*$, where $M= (\bold F_{k|l}
(\pi\mathfrak{sp}_{n|n}(\mathbb{C})))_{\red}$. In particular,
$$
H^0(\bold F_{k|l} (\pi\mathfrak{sp}_{n|n}(\mathbb{C})))\simeq
\bigwedge V_1^*\simeq \bigwedge ((n+1)n/2).
$$

Assume that $\gamma|_{V_2}$ is injective. Then $a_i+a_j< 0$ for all
$i< j$. Hence, $k|l=(0,\ldots, 0|n,l_2,\ldots, l_r)$ or $(0,\ldots,
0|n-1,l_2,\ldots, l_r)$ or $(1,0,\ldots, 0|n-1,n-1,l_3,\ldots,
l_r)$. In these cases
$(\mathfrak{p}(x)_{\pi\mathfrak{sp}_{n|n}(\mathbb{C})})_{\bar
1}\subset V_1$ and $V^*_2$ is the maximal $\mathfrak{g}_{\bar
0}$-module in $\Im \gamma^*$. By Proposition \ref{funk_dop} it
follows that
 $H^0(\bold F_{k|l} (\pi\mathfrak{sp}_{n|n}(\mathbb{C})))\simeq
\bigwedge V_2^*$.

If  $\gamma|_{V_k}$, $k=1,2$, is not injective, then $H^0(\bold
F_{k|l} (\pi\mathfrak{sp}_{n|n}(\mathbb{C})))\simeq \mathbb{C}$ by
Lemma~\ref{gamma is inecive} and Theorem \ref{funk1}.

Suppose that $n=1$. Then
$$
\pi\mathfrak{sp}_{1|1}(\mathbb{C}) = \left\{
\left(
  \begin{array}{cc}
    x & 0 \\
    z & -x \\
  \end{array}
\right), \,\,x,z\in \mathbb{C} \right\}
$$
Since $k_1+l_1\leq n=1$, we have $k|l=(1|0)$ or $(0|1)$. In the
first case, $\mathfrak{p}(x)_{\pi\mathfrak{sp}_{1|1}(\mathbb{C})} =
\operatorname{diag}(x,-x)$ and $\bold F_{1|0}
(\pi\mathfrak{sp}_{1|1}(\mathbb{C}))\simeq (\pt, \bigwedge (1))$. In
particular, $ H^0(\bold F_{1|0}
(\pi\mathfrak{sp}_{1|1}(\mathbb{C}))) \simeq \mathbb{C} \oplus
\mathbb{C}$.

In the second case, $
\mathfrak{p}(x)_{\pi\mathfrak{sp}_{1|1}(\mathbb{C})} =
\pi\mathfrak{sp}_{1|1}(\mathbb{C}))$ and $\bold F_{0|1}
(\pi\mathfrak{sp}_{1|1}(\mathbb{C}))\simeq (\pt,\mathbb{C})$ In
particular, $H^0(\bold F_{0|1} (\pi\mathfrak{sp}_{1|1}(\mathbb{C})))
\simeq \mathbb{C}$.$\Box$

\medskip

\noindent{\bf Case $\mathfrak{g}=\mathfrak{q}_{n|n}(\mathbb{C})$.}
Let $e_1,\ldots,e_n,\pi(e_1),\ldots, \pi(e_n)$ be a basis of
$\mathbb{C}^{n|n}$ which agrees with $\pi$. Consider the superdomain
$Z_{I}$ in $\bold F_{k|k}(\mathfrak{q}_{n|n}(\mathbb{C}))$
corresponding to $I_{s\bar 0}= I_{s\bar 1} = (1,\ldots,k_s)$. Denote
by $x$ the origin of $Z_{I}$. We see that $x=(V_1,\ldots,V_r)$,
where $V_i=\langle e_1,\ldots,e_{k_i}\rangle \oplus \langle
\pi(e_1),\ldots, \pi(e_{k_i})\rangle$. Denote by
$\mathfrak{p}(x)_{\mathfrak{q}}$ the Lie superalgebra of the
stabilizer of $x$ for the action (\ref{dey}) of
$\Q_{n|n}(\mathbb{C})$. It is easy to see that
$$
\mathfrak{p}(x)_{\mathfrak{q}} =\{X\in
\mathfrak{q}_{n|n}(\mathbb{C})\,\mid \, X(V_i)\subset V_i\}.
$$

Let us take a Cartan subalgebra $\mathfrak{t}\subset
\mathfrak{q}_{n|n}(\mathbb{C})_{\bar 0}$ of the following form
$$
\begin{array}{l}
\operatorname{diag}(x_1,\ldots,x_n,x_1,\ldots,x_n).
\end{array}
$$
The corresponding root system $\Delta=\Delta_{\bar 0} \cup
\Delta_{\bar 1}$ is given by
$$
\begin{array}{l}
\Delta_{\bar 0}=\{x_i- x_j\, \, \mid\, i\ne j \},\,\,\, \Delta_{\bar
1}= \{x_i- x_j, \,\, \mid\, i\ne j\}.
\end{array}
$$
Let us take an $n$-tuple $a=(a_1,\ldots,a_n)$ of real numbers such
that
$$
a_1=\cdots=a_{k_r}> \cdots>
a_{k_{2}+1}=\cdots=a_{k_1}>a_{k_{1}+1}=\cdots=a_{n}.
$$
Then $(a,a)\in \mathfrak{t}$. Let
\begin{equation}
\label{parabolic_q}
\mathfrak{p}(a)_{\mathfrak{q}}=
\mathfrak{t}\oplus \bigoplus_{\alpha\in \Delta,\,\alpha(a,a)\geq 0}
\mathfrak{q}_{n|n}(\mathbb{C})_{\alpha}.
\end{equation}
Again $\mathfrak{p}(a)$ depends only on the numbers $k_i$,
$i=1,\ldots,r$. From \cite{Ivan_Onishch}, Chapter $4$, $\S\,4$,
Theorem $4.4$, it can be deduced that
$\mathfrak{p}(a)_{\mathfrak{q}} = \mathfrak{p}(x)_{\mathfrak{q}}$.


\medskip

\t\label{theor_funk gl} {\it $H^0(\bold
F_{k|k}(\mathfrak{q}_{n|n}(\mathbb{C})))\simeq \mathbb{C}$.

}

\medskip

\noindent{\it Proof.} Since $V:=\mathfrak{q}_{n|n}(\mathbb{C})_{\bar
1}$ is an irreducible $\mathfrak{q}_{n|n}(\mathbb{C})_{\bar
0}$-module and $\mathfrak{p}(a)_{\mathfrak{q}}\cap
\mathfrak{q}_{n|n}(\mathbb{C})_{\bar 1}\ne \{0\}$ for all $a$, the
map $\gamma|_V$ cannot be injective. Now our assertion follows from
Lemma \ref{gamma is inecive} and Theorem \ref{funk1}.$\Box$

\medskip

\noindent  {\bf Acknowledgment.} The author is grateful to A.L.
Onishchik, P. Heinzner, A. T. Huckleberry for useful discussion.

\smallskip

\noindent {\textsc{Ruhr-Universit\"{a}t Bochum,
Universit\"{a}tsstra{\ss}e 150, 44801 Bochum, Germany;}}


 \noindent {\emph{E-mail address:}
\verb"VishnyakovaE@googlemail.com"}

\end{document}